\documentclass {article}
\usepackage{graphicx}
\usepackage{amssymb,amsthm,amsmath,harvard}

\catcode`@=11 \catcode`!=11

\expandafter\ifx\csname fiverm\endcsname\relax
  \let\fiverm\fivrm
\fi
  
\let\!latexendpicture=\endpicture 
\let\!latexframe=\frame
\let\!latexlinethickness=\linethickness
\let\!latexmultiput=\multiput
\let\!latexput=\put
 
\def\@picture(#1,#2)(#3,#4){%
  \@picht #2\unitlength
  \setbox\@picbox\hbox to #1\unitlength\bgroup 
  \let\endpicture=\!latexendpicture
  \let\frame=\!latexframe
  \let\linethickness=\!latexlinethickness
  \let\multiput=\!latexmultiput
  \let\put=\!latexput
  \hskip -#3\unitlength \lower #4\unitlength \hbox\bgroup}

\catcode`@=12 \catcode`!=12

\input pictex.tex

\catcode`@=11 \catcode`!=11
  
\let\!pictexendpicture=\endpicture 
\let\!pictexframe=\frame
\let\!pictexlinethickness=\linethickness
\let\!pictexmultiput=\multiput
\let\!pictexput=\put

\def\beginpicture{%
  \setbox\!picbox=\hbox\bgroup%
  \let\endpicture=\!pictexendpicture
  \let\frame=\!pictexframe
  \let\linethickness=\!pictexlinethickness
  \let\multiput=\!pictexmultiput
  \let\put=\!pictexput
  \let\input=\@@input   
  \!xleft=\maxdimen  
  \!xright=-\maxdimen
  \!ybot=\maxdimen
  \!ytop=-\maxdimen}

\let\frame=\!latexframe

\let\pictexframe=\!pictexframe

\let\linethickness=\!latexlinethickness
\let\pictexlinethickness=\!pictexlinethickness

\let\\=\@normalcr
\catcode`@=12 \catcode`!=12

\citationstyle{agsm}

\title{Fictitious Play in $3\times 3$ Games: \\ the transition between periodic and chaotic  behaviour}

\author{Colin Sparrow,
Sebastian van Strien and 
Christopher Harris}

\def\text#1{\mbox{#1}}

\newtheorem {lemma}{Lemma}[section]
\newtheorem {theo}{Theorem}[section]

\newtheorem {prop}{Proposition}[section]

\newtheorem {coro}{Corollary}[section]

\renewcommand{\rho}{\varrho}
\newcommand\st{\,\,;\,\,}

\newcommand\rz{{\mathbb R}}
\newcommand\PP{{\mathbb P}}

\begin{document}

\maketitle

\centerline{Mathematics Institute, University of Warwick and 
Department of Economics, Cambridge}
\centerline{}
\centerline{csparrow@maths.warwick.ac.uk, strien@maths.warwick.ac.uk and 
cjh22@cus.cam.ac.uk}

\begin{abstract}
{In the 60's Shapley provided an example of a
two player fictitious game with  periodic behaviour.
In this game,  player $A$ aims to copy $B$'s behaviour
and player $B$ aims to play one ahead  of player $A$.
In this paper we generalize Shapley's example by introducing
an external parameter. We show that the periodic behaviour
in Shapley's example at some critical parameter value disintegrates
into  unpredictable (chaotic) behaviour, with players
dithering a huge number of times between different strategies. 
At a further critical parameter 
the dynamics becomes periodic again, but now both players
aim to play one ahead of the other.
In this paper we adopt a geometric (dynamical systems) approach
and show that GENERICALLY the dynamics of fictitious play
is essentially continuous and uniquely defined everywhere except at the interior equilibrium. 
Here we concentrate on
the periodic behaviour, while in the sequel to this paper we shall describe the  chaotic behaviour.
} 
\end{abstract}




\section{Introduction}

Suppose that a two-player finite game in normal form is repeated infinitely
often. Suppose further that, at each stage of the repeated game, each player
plays a best response to the empirical distribution of the past play of her
opponent. Then the state of the repeated game at any given stage can be
summarized by the pair consisting of the empirical distribution of past
actions of player $A$ and the empirical distribution of past actions of
player $B$, and this state evolves according to a dynamical system. This
dynamical system is known as fictitious play.

Fictitious play arises in many different guises in economics. It was
originally proposed as an algorithm for finding equilibrium, see  
\citeasnoun{Brown49}, \citeasnoun{Brown51a} and 
\citeasnoun{Brown51b}. It was later
reinterpreted as a model of boundedly rational learning in interactive
settings, 
see Chapter 2 of 
\citeasnoun{Fudenberg-Levine98} and the
references cited there. Several modifications have been introduced,
including: (i) a modification in which players do not play exact best
responses but instead choose actions with a probability that increases with
their payoffs, see \citeasnoun{Cowan92} and 
\citeasnoun{Fudenberg-Levine95};
(ii) a modification in which each player is replaced by a population of
players whose payoffs differ by a small random perturbations, 
see \citeasnoun{Ellison-Fudenberg00}; and (iii) a
modification in which players play exact best responses with some
probability and experiment at random with some probability, see 
\citeasnoun{Fudenberg-Levine98} and the references cited
there. Finally, there is a discrete-time version, in which the underlying
normal-form game is played for each $t\in \left\{ 0,1,2,...\right\} $, and a
continuous-time version, in which the underlying normal-form game is played
for each $t\in \left[ 0,+\infty \right) $.

Much effort has been devoted to finding sufficient conditions under which
fictitious play converges to an equilibrium from all starting points. For
example, fictitious play has been shown to converge in: (i) zero-sum games, see  \citeasnoun{Robinson51}, 
\citeasnoun{Hofbauer95}  and 
\citeasnoun{Harris98};  (ii) generic 
$2\times n$ nonzero-sum games (see 
\citeasnoun{Miyasawa61} and 
\citeasnoun{Metrick-Polak94}
for the $2\times 2$ case; 
\citeasnoun{Sela00} for the $2\times 3$ case; and
\citeasnoun{Berger2005} for the general case); (iii) potential games
(see 
\citeasnoun{Monderer-Shapley96} for the case of weighted
potential games and  
\citeasnoun{Berger2003} for the case of ordinal
potential games); (iv) dominance-solvable games (see 
\citeasnoun{Milgrom-Roberts91}) and (v) various classes of
supermodular games (see: 
\citeasnoun{Krishna92} for the case of supermodular games with
diminishing returns; 
\citeasnoun{Hahn99} for the case of $3\times 3$
supermodular games; and 
\citeasnoun{Berger2003} for the case of
quasi-supermodular games with diminishing returns, the case of
quasi-supermodular $3\times n$ games and the case of quasi-supermodular $%
4\times 4$ games).

Some effort has also been devoted to understanding fictitious play in games
in which fictitious play does not necessarily converge to an equilibrium.
One approach is by way of examples.
The first contribution to this strand of the literature was made by 
\citeasnoun{Shapley64}, who identified a generic $3\times 3$ nonzero-sum game for
which fictitious play has a stable limit cycle. Three more recent
contributions are:  \citeasnoun{Jordan93}, who constructed a 
$2\times
2\times 2$ game for which fictitious play has a stable limit cycle;
\citeasnoun{GH95}, who established an interesting
relationship between the replicator dynamics and fictitious play in three
examples in which fictitious play converges to a limit cycle (namely
Shapley's example, Jordan's example and the rock-scissors-paper game); and
\citeasnoun{Cowan92}, who identified a difficult numerical example in which
fictitious play exhibits chaotic behaviour (in the sense that the flow
contains a subshift of finite type). 
Our work is closest in spirit to  
\citeasnoun{Berger1995}.\footnote{We would like to thank the referee for bringing this unpublished work to our attention. Berger considers a  family of symmetric bimatrix games
depending on a parameter $k$ such that
(i) for $k\in (-2,0)$ has a Shapley orbit  which
degenerates as $k\to 0$, and (ii) for $k>0$, as in our case, there exists 
a hexagonal orbit along which both players are indifferent between (at least) 2 strategies. Berger observes similar 'chaotic' numerical phenomena as in  Subsection~\ref{subsec:chaos} of our paper, and also shows the existence of an additional periodic orbit of saddle-type.}
  
Another approach is taken by \citeasnoun{Krishna-Sjostrom98}, 
building on earlier work of \citeasnoun{Rosenmuller71}.
In this interesting paper, they show that  generically (i.e. for Lebesgue almost all pay-off matrices), fictitious play cannot converge cyclically to a mixed strategy
equilibrium if both players use more than two pure strategies.
In other words, if fictitious play converges to a mixed strategy
equilibrium with both players using more than two strategies 
then the choice of strategies cannot follow a cyclic 
pattern.\footnote{For their result Krishna and Sj\"ostr\"om only (need to) consider the situation where at any given moment only one of the players is in indifferent. As becomes clear from our paper (and its sequel), it is precisely near the set where {\em both} players are simultaneously indifferent that much of the interesting
behaviour happens (and this set 'organises' the local dynamics).}

The purpose of this 
and a sequel to this paper is a detailed analysis of the following family of  $3\times 3$ games with utilities
$v^A$ and $v^B$
determined by the matrices 
\begin{equation}
A= \left(\begin{array}{ccc} 1 & 0 & \beta \\ \beta & 1 & 0 \\ 0 & \beta & 1
\end{array} \right) \quad
B=\left(\begin{array}{ccc}  -\beta & 1 & 0 \\ 0 & -\beta & 1 \\ 1 & 0 & -\beta
\end{array} \right),
\end{equation}
which depend
on a parameter $\beta\in (-1,1)$ (and best response 
dynamics given by the differential inclusion (\ref{eq1})).
We will show that for some parameters $\beta$
play becomes cyclic and for others it becomes erratic.
%

For $\beta=0$, this game
is equal to Shapley's: irrespective of the starting position,
every orbit tends to a periodic motion, 
with eventually player $A$ copying
the previous move of player $B$, while player $B$ is aiming to play 
one ahead  from the previous move of player $A$.
Plotting the index of the 
largest component of the utilities $v^A$ and $v^B$
(``the strategy'' of player $A$ resp. $B$) 
as a function of time we obtain the  repeating pattern (repeating in 6
steps):
$$ 
\begin{array}{cccccccc}	
\mbox{ time duration \quad } & 0.33123 & 0.31767 & 0.33123 & 0.31767 & 0.33123 & 0.31767  &  \\
\mbox{strategy player A \quad }  &1 & 2 & 2 & 3  & 3 & 1 & \\ 
\mbox{strategy player B \quad }  & 2 & 2 & 3 & 3 & 1 & 1 & 
\end{array}
$$
That is, during 0.33123 time units players $A,B$ prefer strategies
$1,2$ respectively, then for 0.31767 time units they prefer
$2,2$ and so on.

For $\beta=\sigma:=(\sqrt{5}-1)/2\approx 0.618$ (ignoring usual
convention we shall call this number\footnote{$\sigma$ is the smaller
of the two roots of $\sigma^2+\sigma-1=0$.}
the golden mean),
the game is equivalent to a zero-sum game\footnote{We would like to 
thank the referee for pointing this out.}:
rescaling $B$ to $\tilde B=\sigma(B-1)$ 
gives $A+\tilde B=0$.
So all orbits tend to the equilibrium $E=(E^A,E^B)$ with
$E^A=(A^B)^T=(1/3,1/3,1/3)$. However, the players will not
choose strategies cyclically (this follows of course 
from  \citeasnoun{Krishna-Sjostrom98}):  one gets sequences such as 
$$\begin{array}{l}	
\mbox{strategy player A \quad }  
         1     2     2     3     3     3     1     1
     2     2     3     3     2     1     1     2
     2     1     1     2     2     3     3     2
     2     3     3     3     2     2     3     3
     1     1     2     2     1     1     3     2
     2     3     3     2     2     3     3     3
     2     2     3     3     1     1     1     2     2     3     3     1     1 \\ 
\mbox{strategy player B \quad }  
    2     2     3     3     2     1     1     2
     2     3     3     1     1     1     2     2
     3     3     2     2     3     3     2     2
     3     3     2     1     1     3     3     1
     1     2     2     3     3     2     2     2
     3     3     2     2     3     3     2     1
     1     3     3     1     1     3     2     2
     3     3     1     1     2\end{array}
$$
Indeed, as will be shown in \citeasnoun{CS2006}, the convergence happens in a chaotic fashion
and can be projected into one which is 'area preserving'.

For $\beta\approx 1$ the players' game also tends
to a periodic pattern,
but now both player $A$ and $B$
aim to play one ahead from the other player,
and therefore now infinitely often repeat  (again repeating in 
6 steps):
$$ 
\begin{array}{cccccccc}	
\mbox{ time duration \quad }  & 0.12060 & 0.39493 & 0.12060 & 0.39493 & 0.12060 & 0.39492  &  \\
\mbox{strategy player A \quad }  & 1 & 1 & 3 & 3  & 2 & 2 & \\ 
\mbox{strategy player B \quad }  & 3 & 2 & 2 & 1 & 1 & 3 & 
\end{array}
$$
We call this the {\em anti-Shapley} pattern.
The main point of this paper is to show how
this change of behaviour occurs. It turns out that 
there are two {\em critical} parameters, $\sigma$ the golden mean 
$\approx 0.618$ and $\tau\approx 0.915$ such that 
\begin{itemize}
\item for  $\beta\in (-1,\sigma)$ the players typically end up repeating
Shapley's pattern;
in fact, for $\beta\in (-1,0]$ regardless of the initial positions,
the behaviour  tends to a periodic one
(for $\beta>0$ this is not true: then there are orbits
which tend to the interior equilibrium and/or other
periodic orbits).
\item for $\beta\in (\sigma,\tau)$ the players
become extremely indecisive and erratic (and the moves
become chaotic), while 
\item for 
$\beta\in (\tau,1)$ the players typically end up playing
the anti-Shapley pattern.
\end{itemize}
For example, when $\beta=0.75$ one sees erratic behaviour such as
\begin{equation*}
\begin{array}{l}
\mbox{strat. A } \quad  3     3     1     1     2     2     2     3     3     1     1     2     2     1     1     3     2     2     3     3     3     2  2     3     3     2     2     1     3     3     1     1     3     3     3     1     1     3     3     1     1     3     3     1  1     3     3     1     1     3     3     2     1     1     2     2     1     1     2     2
\\
\mbox{strat. B }  \quad 2     1     1     2     2     1     3     3     1     1     2     2     3     3     2     2     2     3     3     2     1     1 3     3     1     1     3     3     3     1     1     3     3     2     1     1     2     2     1     1     2     2     1     1
2     2     1     1     2     2     1     1     1     2     2     1     1     2     2     1
\end{array}\end{equation*}
In this paper we shall describe the bifurcations that the periodic
orbits undergo as $\beta$ varies; in a sequel to this paper
(\citeasnoun{CS2006}),
the nature of the chaotic motion will be described in 
more detail (for example,  that even for $\beta\in (0,\sigma)$, there are many (exceptional) orbits that are attracted to the equilibrium point).

Of course, when one or more of the players is 
indifferent between two strategies
their dynamics is not uniquely determined. However, as we shall
show in Corollary~\ref{cor:uniqflow} (in
Subsection~\ref{ss:cd1}) under a mild assumption
{\em only exceptional} orbits of the play are affected by this ambiguity.
In particular, the orbits discussed above are 
uniquely determined.

\section{Continuous Time Fictitious Play: the Generic Case}
\label{ss:ctfp}
We study the following version of continuous time
fictitious play, which is convenient for calculation, for
numerical simulation, and for extracting the geometrical properties of
the dynamics. It is equivalent to other formulations in the literature.
Two players A and B, both having
$n$ strategies to choose from, have $n \times n$ pay-off 
matrices $A = (a_{ij})$ and $B = (b_{ij})$, $1\le i,j \le n$ respectively.
(The context will always make it clear whether we are talking about Players
A or B or about the matrices $A$ and $B$).
We will assume throughout this paper
that the matrices $A$ and $B$ are both {\em non-singular} (i.e. have
non-zero determinant).
We denote the set of probability vectors for players A and B by
$\Sigma_A, \ \Sigma_B \subset \rz^n$. (Though both simplices are the same
subset of $\rz^n$, the subscripts are useful when we wish
to make it clear which players probabilities we are discussing).
At some time  $t_0 \in \rz$
the {\it state} of the game is given by a pair of probability
vectors, $(p^A(t_0),\ p^B(t_0)) \in 
\Sigma:= \Sigma_A \times\Sigma_B$, representing 
the current (usually mixed) strategies of the two players.  It is convenient
to define $p^A$ and $p^B$ so that they are row and column vectors 
respectively.

The dynamics
of the system is determined as follows.   
At time $t_0$,  each player computes a
{\it best response} $BR_A(p^B) \in \Sigma_A$
and $BR_B(p^A) \in \Sigma_B$ to the other player's 
current strategy.  If we denote the set of pure strategies for the
two players by
$P^A_i \in \Sigma_A$ and $P^B_j \in \Sigma_B$, $1 \le i,j,\le n$,
then  generically
\footnote{In this paper the word `generic' means that one is not in a codimension-one
subspace}
player A's  best response $BR_A$ is a pure strategy $P^A_i\in \Sigma_A$,
for some 
$i \in \{1,\dots,n\}$
chosen so that the $i$-th coordinate of the vector 
$$v^A(t_0) = Ap^B(t_0)$$  
is larger than the other coordinates (the non-generic case, when
several coordinates of the vector $v^A(t_0)$ are equal will be treated later on). Similarly,
player B's best response $BR_B$ is generically a pure strategy
$P^B_j\in \Sigma_B$, $j \in \{1,\dots,n\}$
where $j$ is chosen so that the $j$-th coordinate of the vector 
$$v^B(t_0)= p^A(t_0)B$$  
is larger than any of the other coordinates. In cases where players are 
indifferent between two or
more strategies (the non-generic case), any convex combination of the relevant pure 
strategies is an equally good best
response, so 
$$
BR_A\colon \Sigma_B\to \Sigma_A\mbox{   and }
BR_B\colon \Sigma_A\to \Sigma_B$$
are in fact {\em set-valued maps}.
This means equation (\ref{eq1}) below actually
defines a differential inclusion rather than a differential equation. 
\footnote{We will deal with this apparent complication 
in sections \ref{ss:cd1} and \ref{ss:mi} below. Our approach will
be to consider the flow where best responses {\it are} uniquely defined
and to investigate the extent to which this flow extends in a unique
and continuous way onto the parts of phase space where best responses 
are not unique; for the moment readers unfamiliar with the idea 
of differential inclusion will lose very little by continuing to think 
of (\ref{eq1}) as a differential equation.}

Once best responses are selected, the dynamics is determined by the 
piecewise linear equations:
\begin{equation}
\begin{array}{ll}
       {dp^A}/{dt} & =  BR_A(p^B) - p^A \\
       {dp^B}/{dt} & =  BR_B(p^A) - p^B \end{array}
\label{eq1}\end{equation}
so that each player's tendency is to adjust their strategy in a 
straight line from their current strategy towards their current 
best response.\footnote{The best response correspondences 
$p^B\mapsto BR_A(p^B)$ and $p^A\mapsto BR_B(p^A)$
are upper semicontinous with values closed, convex sets.
It follows from 
\citeasnoun[Chapter 2.1]{AubinCellina84}
that through each initial value there exists at least one solution
which is Lipschitz continuous and defined for all positive 
time. In a sequel to this paper  we will show that, under mild regularity 
conditions, any solution is in fact piecewise linear.}
In the GENERIC case as discussed above, player A continuously
adjusts her probability vector $p^A(t)$ from its current position $p^A(t_0)$ 
in a straight line towards $P^A_i$, and player B similarly adjusts $p^B(t)$ to
move in a straight line towards $P^B_j$, so that for times greater than 
$t_0$, and for so long as strategies $P^A_i$ and $P^B_j$ 
remain unique best responses,  
we may write, after a reparameterisation of time, the solution of 
(\ref{eq1}) becomes
\begin{equation}
\begin{array}{ll}
p^A(t_0+s)&=p^A(t_0)(1-s) + s \cdot P^A_i,\\
p^B(t_0+s)&=p^B(t_0)(1-s) +  s \cdot P^B_j,
\end{array}
\label{eq2}
\end{equation}
for $s \in [0,1]$.  
Note that we have parameterised time 
so that both players move towards their best response strategies
at a uniform speed, such that if their best responses did
not change they would arrive at $(P^A_i, P^B_j)\in \partial \Sigma$ at time $t=t_0+1$, 
at which point
the dynamics would halt.  The parameterisation is
chosen for ease of calculation and computation, and does not affect the
geometry of the dynamics which is our chief concern in this paper; it would
be straightforward to use the original time $\rho = -\ln (1-s)$ if
an exponential approach to $(P^A_i, P^B_j)$ in infinite time were preferred
(the speed of a solution of (\ref{eq1}) is proportional to the distance to the target point).  
In any case, since orbits consist of pieces of line segments, numerical simulations
can be done easily and with arbitrary accuracy.

Equations (\ref{eq2}) therefore determine the dynamics up {\em until such time}
as one or other (or both) players become {\em indifferent} 
between two (or more) pure strategies. 

Let us first describe the sets where
the best responses are constant. We will be particularly interested in games
where there is an {\it interior equilibrium}; 
that is a point $E=(E^A, E^B)$ in the
interior of $\Sigma_A \times \Sigma_B$ where both players are 
indifferent between {\em all} $n$ strategies, so that all components
of $AE^B$ are equal, and all components of $E^AB$ are equal.
The following lemma is entirely straightforward.

\begin{lemma}[The complement of the indifferent sets]
 \label{convex}
Assume that $A$ and $B$ are $n\times n$ non-singular matrices.
Then the set where one player is indifferent between two given strategies (and so where she has a choice of best responses)
forms a codimension-one hyperplane and for each $1\le i,j\le n$
the set  $$S_{ij}= \{(p^A,p^B)\st BR_A(p^B)=P^A_i \mbox{ and } 
BR_B(p^A)=P^B_j\}$$
is convex.

If there exists an equilibrium  $E=(E^A, E^B)$
in the interior of $\Sigma$, 
then all $n^2$ regions $S_{ij}$ are non-empty, and their closures 
meet at $(E^A, E^B)$.
\end{lemma}

\bigskip

\begin{figure}[htp]
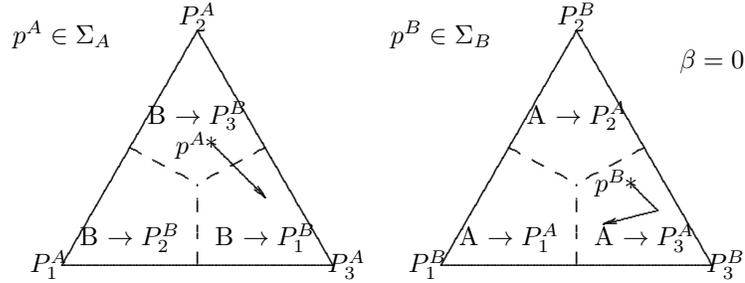
 \hfil
\beginpicture
\dimen0=0.18cm
\setcoordinatesystem units <\dimen0,\dimen0>
\setplotarea x from 0 to 20, y from -5 to 20
\setlinear
\plot 0 0 20 0 /
\plot 0 0 10 17.3 20 0 /
\put {$p^A$} at 9.5 8.8
\put{$*$} at 11 9
\arrow <5pt>  [0.2,0.4] from 11 9 to  15 5
\setdashes
\plot 10 0 10 6 /
\plot 5 8.66 10 6 /
\plot 15 8.66 10 6 /
\put {$P^A_1$} at -1 0
\put {$P^A_2$} at 10 18.3
\put {$P^A_3$} at 21 0
\put {B $\to P^B_2$} at 5 2
\put {B $\to P^B_1$} at 15 2
\put {B $\to P^B_3$} at 10 11
\put {$p^A \in \Sigma_A$} at 0 17
\setcoordinatesystem units <\dimen0,\dimen0> point at -28 0
\setplotarea x from 0 to 30, y from -5 to 20
\setsolid
\put {$p^B \in \Sigma_B$} at 0 17 \put {$\beta=0$} at 20 15
\setlinear
\plot 0 0 20 0 /
\plot 0 0 10 17.3 20 0 /
\put {$p^B$} at 12.5 6
\put {$*$} at 14 6
\plot 14 6 16 4 /
\arrow <5pt>  [0.2,0.4] from 16 4 to  12 3 
\setdashes
\plot 10 0 10 6 /
\plot 5 8.66 10 6 /
\plot 15 8.66 10 6 /
\put {$P^B_1$} at -1 0
\put {$P^B_2$} at 10 18.3
\put {$P^B_3$} at 21 0
\put {A $\to P^A_1$} at 5 2
\put {A $\to P^A_2$} at 10 11
\put {A $\to P^A_3$} at 15 2
\endpicture
\caption{\label{fig:ssb=0}
{\small An example in which the simplices 
$\Sigma_A$ and $\Sigma_B$
the regions where players have a specific preference
are as shown, separated by dashed lines which show the 
indifference sets $Z^B \subset \Sigma_A$ where player B is indifferent
between two or more strategies, and $Z^A \subset \Sigma_B$ where player A is
indifferent between two or more strategies. 
Starting from an initial condition $(p^A,p^B)$ (marked with the * symbols), player A moves
towards $P^A_3$, while Player B starts by moving towards $P^B_3$ 
and changes direction towards $P^B_1$ as $p^A$ crosses $Z^B$. 
The sets $S_{ij}\in \rz^4$ are the sets formed as products of two
kite-shaped regions where each player has a strict preference.}
}
\end{figure}

Figure~\ref{fig:ssb=0} shows a $3\times 3$ example in which
the 
two
simplices $\Sigma_A$ and $\Sigma_B$ are each two dimensional triangles,
with vertices representing the pure
strategies of the two players as labelled, and the regions $S_{ij}$
with different aims marked in the figure. The state of the
system at a particular time, $(p^A(t), p^B(t))$ is represented by
a pair of points, $p^A$ and $p^B$, one in each simplex. The local 
direction of motion is towards the vertex determined by the
location of the other players position in her simplex, as
illustrated. Motion will continue in a straight line until
one or other player hits an indifference line (see caption). 

\subsection{Uniqueness and continuity of dynamics where at most one 
player is indifferent}
\label{ss:cd1}

We now consider, for general matrices $A$ and $B$, the hyperplanes 
where one or other player is 
indifferent between two strategies.

Let us denote the set where player A is indifferent
between strategies $P^A_i$ and $P^A_j$ 
by $Z^A_{i,j} \subset \Sigma_B$. We
define $Z^B_{i,j} \subset \Sigma_A$ similarly. If we set
$Z^A =  \cup_{i,j} Z^A_{i,j}$, then 
$\Sigma_A \times Z^A$ is the set
where player A does not have a unique best response. We 
define the set $Z^B$ similarly. 
Outside the set 
$$Z=\left( \Sigma_A\times Z^A\right) \bigcup
\left( Z^B \times \Sigma_B \right)  $$ 
where at least one of the players is indifferent between two or more
strategies the flow is uniquely
defined by (\ref{eq1}) and is clearly continuous. 
If we define 
\begin{equation}
Z^*=Z^B \times Z^A \subset \Sigma_A\times \Sigma_B\ ,
\end{equation}
which is the subset of $Z$ where {\it both} players are indifferent between at 
least two strategies, then the next lemma says that 
the uniqueness and continuity extends to the {\em complement}
of $Z^*$ provided matrices $A$ and $B$ satisfy 
the condition (\ref{eqn:trans}) in the lemma.
(Note that this transversality condition is satisfied for
the Shapley example from above if and only if  $\beta\ne 0$.)

\medskip
\begin{prop} [Transversality condition implies local continuity and uniqueness outside
the set where both players are indifferent] \label{contitra}
The motion defined by equations (\ref{eq1})
forms a continuous flow on  
$(\Sigma_A\times \Sigma_B)\setminus Z^*$
provided the following condition is met:
for any point $(p^A,p^B)\notin Z^*$ with, say,
$p^A \in Z^B$, $p^B \notin Z^A$ and 
strategy $P^A_k$ preferable for player A,
\begin{equation}
\mbox{ the vector $P^A_k$ is not parallel to the
plane $Z^B\subset \Sigma_A$ at $p^A$}
\label{eqn:trans}
\end{equation}
(and similarly, in the case when the role of $p^A$ and $p^B$ is reversed).
\footnote{The linear space of vectors where the $l$ and $l'$-th
options are indifferent to player $B$ consists of the
vectors in ${\mathbb R}^m$ which are orthogonal
to the difference of the $l$-th and $l'$-th column 
vector of the matrix $B$.
That the $k$-unit unit vector in ${\mathbb R}^m$ is not contained
in this, means therefore that $b_{kl}\ne b_{kl'}$.
So the above transversality condition holds if
within each column vector of $A$, the coefficients
are all distinct (i.e. $a_{ij}\ne a_{i'j}$ when $i\ne i'$), and similarly
$b_{kl}\ne b_{kl'}$ when $l\ne l'$.  
Clearly the set of matrices for which this holds is open and dense
(and of full Lebesgue measure).} 
\end{prop}
\begin{proof}
The transversality condition (\ref{eqn:trans}) ensures that a point is
moved immediately off the codimension-one plane
$Z^B\times \Sigma_B$ so
that player B has a unique best response the moment the motion begins.
\end{proof}

Since $Z^*$ has zero Lebesgue measure and 
consists of codimension two planes, this shows

\begin{coro}\label{cor:uniqflow}
Assume the tranversality condition (\ref{eqn:trans}) holds.
Then there exists a set $X\subset \Sigma_A\times \Sigma_B$
of initial conditions such that
\begin{itemize}
\item for each starting point in $X$ the flow
defined by (\ref{eq1}) is continuous and unique  for all $t\ge 0$.
\item $X$ has full Lebesgue measure and is open and dense.
\end{itemize}
\end{coro}

\subsection{Local dynamics where both players are indifferent: outside sets of codimension three}
\label{ss:mi}

Even though the previous corollary shows that
most starting positions never enter the set $Z^*$
where both players are indifferent, it is
still important to analyse the dynamics {\em near} this set.
As it turns out, in many places the flow outside $Z^*$
extends continuously to $Z^*$. This suggests
a  natural choice of the non-uniquely defined dynamics within $Z^*$.
Much more importantly - even if we are not interested
in orbits that lie on $Z^*$ and do not want to make a choice
for the dynamics on this set -  
the continuous extension dynamics on $Z^*$ is an
``organising principle'' for nearby orbits, as we shall
see below. This continuous extension dynamics
makes that all orbits near $Z^*$ have consistent
dynamics.  This point of view puts us
in a classical dynamical systems framework and  allows us, for example,  to describe bifurcations of the flow (i.e. changes in
 dynamics
when the parameter $\beta$ changes) for orbits near $Z^*$
(as we shall see in the Section~\ref{s:simpleprops}).

\begin{figure}[htp]
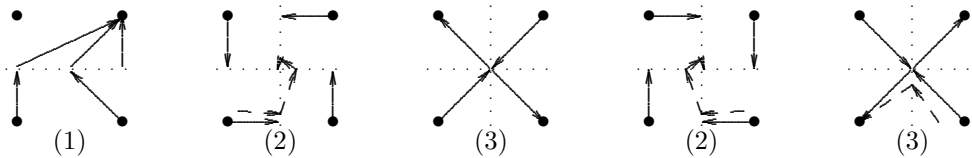
 \hfil
\beginpicture
\dimen0=0.07cm \dimen1=0.07cm
\setcoordinatesystem units <\dimen0,\dimen1> point at 40 0
\setplotarea x from -10 to 10, y from -15 to 15
\setlinear \setdots
\plot -12 0 12 0 /
\put {$(1)$} at 0 -14
\setsolid
\put {$\bullet$} at 10 10
\put {$\bullet$} at -10 10
\put {$\bullet$} at -10 -10
\put {$\bullet$} at 10 -10
\arrow <5pt> [0.2,0.4] from 10 -10 to 0.3 -0.5
\arrow <5pt> [0.2,0.4] from 0.3 0.5 to 9.5 9.5
\arrow <5pt> [0.2,0.4] from -10 -10 to -10 -0.5
\arrow <5pt> [0.2,0.4] from -10 0.5 to 9.5 9.5
\arrow <5pt> [0.2,0.4] from 10 0.5 to 10 10
\setcoordinatesystem units <\dimen0,\dimen0> point at 0 0
\setplotarea x from -10 to 10, y from -10 to 10
\setlinear \setdots
\plot -12 0 12 0 /
\plot 0 12 0 -12 / \setsolid
\put {$(2)$} at 0 -14
\put {$\bullet$} at 10 10
\put {$\bullet$} at -10 10
\put {$\bullet$} at -10 -10
\put {$\bullet$} at 10 -10
\arrow <5pt> [0.2,0.4] from 10 10 to 0.5 10
\arrow <5pt> [0.2,0.4] from -10 10 to -10 0.5
\arrow <5pt> [0.2,0.4] from -10 -10 to -0.5 -10
\arrow <5pt> [0.2,0.4] from 10 -10 to 10 -0.5
\setdashes
\arrow <5pt> [0.2,0.4] from -8 -8 to 0 -8.5
\arrow <5pt> [0.2,0.4] from 0 -8.5 to 3 0
\arrow <5pt> [0.2,0.4] from 3 0 to 0 2
\arrow <5pt> [0.2,0.4] from 0 2 to -0.5 0
\setsolid
\setcoordinatesystem units <\dimen0,\dimen0> point at -40 0
\setplotarea x from -10 to 10, y from -10 to 10
\setlinear \setdots
\plot -12 0 12 0 /
\plot 0 12 0 -12 / \setsolid
\put {$(3)$} at 0 -14
\arrow <5pt> [0.2,0.4] from 9.5 9.5 to 0.3 0.3 
\arrow <5pt> [0.2,0.4] from -9.5 -9.5 to -0.3 -0.3 
\arrow <5pt> [0.2,0.4] from -0.3 0.3 to -10 10 
\arrow <5pt> [0.2,0.4] from 0.3 -0.3 to 10 -10 
\put {$\bullet$} at 10 10
\put {$\bullet$} at -10 10
\put {$\bullet$} at -10 -10
\put {$\bullet$} at 10 -10
\setcoordinatesystem units <\dimen0,\dimen0> point at -80 0
\setplotarea x from -10 to 10, y from -10 to 10
\setlinear \setdots
\plot -12 0 12 0 /
\plot 0 12 0 -12 / \setsolid
\put {$(2)$} at 0 -14
\put {$\bullet$} at 10 10
\put {$\bullet$} at -10 10
\put {$\bullet$} at -10 -10
\put {$\bullet$} at 10 -10
\arrow <5pt> [0.2,0.4] from 10 10 to 10 0.3
\arrow <5pt> [0.2,0.4] from -10 10 to -0.3 10
\arrow <5pt> [0.2,0.4] from -10 -10 to -10 -0.3
\arrow <5pt> [0.2,0.4] from 10 -10 to 0.3 -10
\setdashes
\arrow <5pt> [0.2,0.4] from 8 -8 to 0 -8.5
\arrow <5pt> [0.2,0.4] from 0 -8.5 to -3 0
\arrow <5pt> [0.2,0.4] from -3 0 to 0 2
\arrow <5pt> [0.2,0.4] from 0 2 to 0.5 0
\setsolid
\setcoordinatesystem units <\dimen0,\dimen0> point at -120 0
\setplotarea x from -10 to 10, y from -10 to 10
\setlinear \setdots
\plot -12 0 12 0 /
\plot 0 12 0 -12 / \setsolid
\put {$(3)$} at 0 -14
\arrow <5pt> [0.2,0.4] from 0.3 0.3 to 10 10 
\arrow <5pt> [0.2,0.4] from -0.3 -0.3 to -10 -10 
\arrow <5pt> [0.2,0.4] from -9.5 9.5 to -0.3 0.3 
\arrow <5pt> [0.2,0.4] from 9.5 -9.5 to 0.3 -0.3 
\setdashes
\arrow <5pt> [0.2,0.4] from 5 -10 to 0 -3
\arrow <5pt> [0.2,0.4] from 0 -3 to -10 -10
\setsolid
\put {$\bullet$} at 10 10
\put {$\bullet$} at -10 10
\put {$\bullet$} at -10 -10
\put {$\bullet$} at 10 -10
\endpicture
\caption{\label{fig:2x2}
{\small The motion in various cases for $2\times 2$ games.}}
\end{figure}

Before going into the general case in the next subsection,
let  us consider the codimension-two case.
To do this,  let us first review the dynamics in the $2\times 2$ case
(this of course is well known, but it is useful to
be explicit here as we shall use the classification later on).
In this case the probability simplices
$\Sigma_A$ and $\Sigma_B$ in $\rz^2$ can both be identified
with $[0,1]$ and so $\Sigma_A\times \Sigma_B$ with the
square $[0,1]\times [0,1]$. 
Player $A$ moves left or right, depending on the position $p^B$ of 
player $B$, i.e. depending whether $(p^A,p^B)$
is above or below a certain horizontal line, and
similarly player $B$ moves up or down depending
on whether $(p^A,p^B)$ is to the left or right
of some vertical line.
Drawn in Figure~\ref{fig:2x2} are a few
orbits of the system in various cases and $\bullet$
marks the aims $(0,0)$, $(1,0)$, $(0,1)$ and $(1,1)$. The dotted
lines denote where one player in indifferent to two strategies.  From
left to right cases (1), (2), (3), (2) and (3).
In cases (2) the orbit spirals towards the interior 
point,
while in case (3) the interior point acts as a 
saddle.\footnote{In case (3), the orbit through a point
in the stable manifold  enters the equilibrium point in finite time!
Moreover, through the interior equilibrium point 
$E$ there are an uncountable set of solutions; 
one of these remains in $E$ for all time (at least locally), and
the others stay in $E$ for some (possibly small) 
finite time which can be  chosen freely, and then fall off in the
unstable direction.}

\begin{lemma} [$2\times 2$ case]
Assume that both players have two strategies
and that the corresponding $2\times 2$ matrices  are non-singular.
Then there are three possibilities:
\begin{enumerate}
\item there is no interior equilibrium in $\Sigma_A\times \Sigma_B$
and all orbits tend to a unique singular point on the boundary
(this occurs if $Z^A$ and or $Z^B$ are empty, and is depicted in 
(1) in Figure~\ref{fig:2x2}); 
\item there is an interior equilibrium $E$ which is stable, orbits
spiral towards it and the flow extends continuously to $E$
(this is depicted in (2) in Figure~\ref{fig:2x2});
\item there is an interior equilibrium $E$ which is a saddle point
and the flow cannot be extended in a continuous way to
$E$ and there is genuine non-uniqueness in the flow near $E$
(this is depicted in (3) in Figure~\ref{fig:2x2}).
\end{enumerate}
\label{lemma:2d}
\end{lemma}
As we shall see in the next section, in case (2)
the solution tends to $E$ as $t\to \infty$,
spiralling faster and faster around $E$ in a very precise 
way.\footnote{In this case the orbit through any point
outside the equilibrium never hits the equilibrium point.}  
This means that we can think of the flow
at $E$ as just stationary there; the flow is continuous (and each solution is unique).

We now consider the behaviour of the system when both players
are indifferent between two or more strategies (but now for general $A$ 
and $B$).
This occurs on the set $Z^* = Z^B \times Z^A$ defined above.
Since we are interested in typical trajectories of the system, and
typical trajectories will not actually hit $Z^*$, we are
particularly concerned with the implications of this for trajectories
near to but not on $Z^*$. This contrasts with the case of 
the codimension one sets
considered above; typical trajectories do hit those.
We will see, though, that we can in many circumstances define a unique
flow for points in $Z^*$ as well.

\begin{figure}[htp]
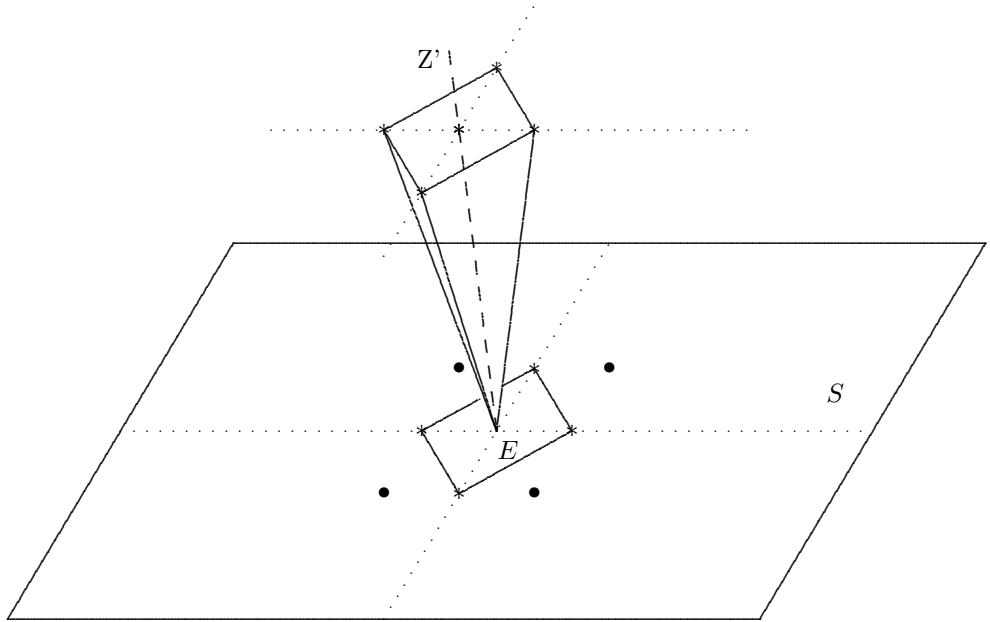
 \hfil
\beginpicture
\dimen0=0.5cm
\setcoordinatesystem units <\dimen0,\dimen0>
\setplotarea x from 0 to 25, y from 0 to 10
\plot 0 0 20 0 26 10 6 10 0 0 /
\put {$S$} at 22 6 
\setdots
\plot 3 5 23 5 /
\plot 10 0 16 10 /
\setsolid
\plot 13 5 11 11.33 14 13 13 5 10 13 11 11.33 /
\plot 10 13 13 14.66 14 13 /
\setdots
\plot 7 13 20 13 /
\plot 10 9.66 14 16.33 /
\multiput {$*$} at 12 13 10 13 11 11.33 14 13 12 13 13 14.66 /
\multiput {$*$} at 12 3.33 15 5 14 6.66 11 5 /
\plot 12 3.33 15 5 14 6.66 11 5 12 3.33 /
\setsolid
\plot 12 3.33 15 5 14 6.66 13.2 6.22 /
\plot 12.5 5.82 11 5 12 3.33 /
\multiput {$\bullet$} at 10 3.33 14 3.33 12 6.66 16 6.66 /
\put {$E$} at 13.3 4.5
\setdashes
\plot 13 5 12 13 11.8 14.6 11.7 15.4 /
\put{Z'} at 11.2 14.9
\endpicture
\caption{\label{fig:quad}
{\small The targets $(P^A_s,P^B_t)$, $s\in \{i,j\}$, $t\in \{k,l\}$ 
in the plane
$S$ are drawn marked as $\bullet$. The dashed line corresponds to
$Z'$. A quadrangle $T_r$ in $S$ is also drawn
and also another quadrangle $T'_r$ which is in
a plane $S_q$ (through $q\in Z'$)
which is parallel to $S$ is drawn.
The cone of the quadrangle over the equilibrium point
$E$ is invariant. Orbits spiral along the cone towards $E$, always
aiming for one of the points marked as $\bullet$. Note that the motion
on the cone is {\bf not} close to the direction of the axis $Z'$. }}
\end{figure}

Consider the subset  $Z'_{kl,ij}$ of $Z^B_{k,l}\times Z^A_{i,j}\subset Z^*$
where  player A is indifferent {\em only} between
strategies $i$ and $j$, and player B is indifferent {\em only} between 
strategies $k$ and $l$.
Since we assumed that $A$ and $B$ are non-singular,
this is a codimension-two plane in $\Sigma_A\times \Sigma_B$.
Define a neighbourhood of $Z'_{kl,ij}$ by
\begin{equation}
N=\left\{(p^A,p^B)\st
BR_A(p^B)\subset \{i,j\} \mbox{ and } BR_{B}(p^A)\subset \{k,l\} \right\}.
\label{eq:setn}
\end{equation}
In the set $N$,
player A prefers strategies $i$ and $j$ to
all other pure strategies, and player B prefers strategies $k$
and $l$. In $N\setminus Z'$, because of Proposition \ref{contitra} above, 
equation (\ref{eq1}) can be justifiably interpreted as a differential 
equation with unique solutions (rather than a differential 
inclusion with possibly non-unique solutions)
with the dynamics consisting of a linear decline
in all components $p^A_m, m\ne i,j$ and
$p^B_{n}, n \ne k, l$ while the behaviour of the
other two components of $p^A$ and $p^B$ depends
entirely on the $2 \times 2$ submatrices of $A$ and $B$:
$$A[k,l]= \left(\begin{array}{cc} a_{k,k} & a_{k,l}  \\ 
                                    a_{l,k} & a_{l,l}  
\end{array} \right) \qquad
B[i,j]= \left(\begin{array}{cc} b_{i,i} & b_{i,j}  \\ 
                                b_{j,i} & b_{j,j}  
\end{array}
\right).\ $$
These matrices describe a game on the two-dimensional space
spanned by $[P^A_i,P^B_j]\times [P^B_{k},P^B_{l}]$
(where $[P^A_i,P^A_j]$ stands for the line segment
$\lambda P^A_i + (1-\lambda)P^A_j\subset \Sigma_A, 0\le \lambda \le 1$).
It follows\footnote{
A similar decomposition is also used in 
\citeasnoun[page 6 and Appendix A]{Hofbauer95}. Amongst other things, Hofbauer shows
in  this interesting paper that the dynamics starting at a point $x\in \Sigma_A\times \Sigma_B$
can proceed in a straight line in the direction $b\in \Sigma_A\times \Sigma_B$ iff $b$ is a Nash equilibrium of the  restricted game at $x$. 
Here we use the same approach to investigate uniqueness.
In a sequel to this paper, we will push this point of view further.}  that the dynamics in $N$
is the product of the
dynamics in this two-dimensional space (on which the dynamics
is as in Figure~\ref{fig:2x2}) and a linear decline along the
codimension-two plane in $Z'$. 
Together with Lemma \ref{contitra} this gives:

\begin{prop}[Continuity and uniqueness outside sets of codimension three]
 \label{prop:contx}
Assume that $A$ and $B$ have maximal rank
and that the transversality condition (\ref{eqn:trans})
is satisfied. Then
the flow defined on the sets $S_{ij}$ has a unique continuous 
extension everywhere
{\em except (possibly) where:}
\begin{itemize}
\item one player is indifferent to at least three strategies
and the other to at least two;
\item both players are indifferent to precisely two strategies
and along this set the dynamics is as in case (3)
described above (in Lemma~\ref{lemma:2d}).
\end{itemize}
In the neighbourhood  $N$ of $Z'_{kl,ij}$ defined in (\ref{eq:setn})  the flow
decomposes into a flow governed by a two by two game, and a constant flow in the remaining direction(s) and there are three cases:
\begin{enumerate}
\item trajectories cross
either $Z^B_{k,l}$ or $Z^A_{i,j}$ and
locally the flow is the product of a constant flow perpendicular
to the plane and a  flow as in (1) in Figure~\ref{fig:2x2};
\item trajectories in $N$ spiral towards $E$ (i.e. towards $Z'$) 
along a cone over a quadrangle
as in Figure~\ref{fig:quad}, and move closer and closer
to the unique trajectory in $Z'$ that remains in $Z'$ -- this trajectory
is therefore the unique continuous extension of the flow in $N$ onto $Z'$ - and locally 
the two-dimensional flow is as in (2) in Figure~\ref{fig:2x2};
\item $Z'$ acts as a saddle and
generic trajectories leave $N$  moving away from $Z'$
in one of two opposite directions and locally 
the two-dimensional flow is as in (3) in Figure~\ref{fig:2x2}.
\end{enumerate}
\end{prop}
\begin{proof}
Only the 2nd case needs some explanation.
Let $S$ be the two dimensional plane spanned by the four targets  
$(P^A_s,P^B_t)$, $s\in \{i,j\}$ and $t\in \{k,l\}$.
Moreover let $V\subset S$ be a quadrangle with
 sides parallel to the lines connecting the equilibrium $E\in S$ 
 to $(P^A_s,P^B_t)$, 
 $s\in \{i,j\}$ and $t\in \{k,l\}$
(in such a way that in the region where the players aim for 
$(P^A_s,P^B_t)$ one takes a line parallel to the line 
connecting $)$ to $(P^A_s,P^B_t)$). Now taking any quadrangle
$V_r$ parallel to $V$, and the cone through the equilibrium 
$E\in S$ over $V_r$. Clearly orbits that start in this cone, 
will stay in it.
\end{proof}


\bigskip




In case (1) though both players have a choice
of best response at $(p^A,p^B)$, there is no choice of best responses
that keeps the trajectory inside $Z^*$. This transversality ensures
that the degeneracy occurs only at discrete times
(i.e., `momentarily' if we may put it like that),
and there is actually only one solution to the differential inclusion 
through $(p^A,p^B)$. In case (2) there is
a unique choice of best responses on $Z^*$ for both players such
that the flow defined by equation (\ref{eq1}) remains within
$Z^*$. Nearby solutions remain close to and
spiral (locally) towards this solution, and so it is the
continuous extension of the (unique) flow off $Z^*$. 
Again, the apparent non-uniqueness
is an artefact; we may say if we like that any attempt to depart from this
solution in $Z^*$ is immediately thwarted by the neighbouring flow. 
Only in case (3) do we
find genuine non-uniqueness of solutions to the differential inclusion. 
This amounts to the fact that there is a choice of best response for 
each player which keeps the flow within $Z^*$, and all other choices 
move you immediately off $Z^*$ in one of two directions. 
Thus through $(p^A,p^B)$ there are an uncountable set of solutions; 
one of these remains in $Z^*$ for all time (at least locally), and
the others stay in $Z^*$ for some (small) finite time which can be 
chosen freely, and then fall off in one of the two directions 
(again chosen freely). 



 

\section{Existence and stability of periodic orbits in
the Shapley family}
\label{s:simpleprops}

From now we shall concentrate on a 
particular family example which includes a classical
example by Shapley. These examples are given by matrices 
\begin{equation}
A= \left(\begin{array}{ccc} 1 & 0 & \beta \\ \beta & 1 & 0 \\ 0 & \beta & 1
\end{array} \right) \quad
B=\left(\begin{array}{ccc}  -\beta & 1 & 0 \\ 0 & -\beta & 1 \\ 1 & 0 & -\beta
\end{array} \right),
\label{eqn:AB}
\end{equation}
When $\beta=0$, this means that the utility of player $A$
is equal to $p^B$, so player $A$'s will aim for
the largest component of $p^B$, while player $B$
utility is $(p^A_3,p^A_2,p^A_1)=p^AB$ so $B$ is aiming
one ahead  from the largest component of $p^A$.
So for $\beta=0$, it is like the paper, scissor, rock 
game. Our interest is to understand how the behaviour
in this game changes when $\beta$ changes.
In Subsection~\ref{subsec:chaos} we show that the
player's actions will become extremely erratic
for certain values of $\beta$.

\begin{figure}[htp]
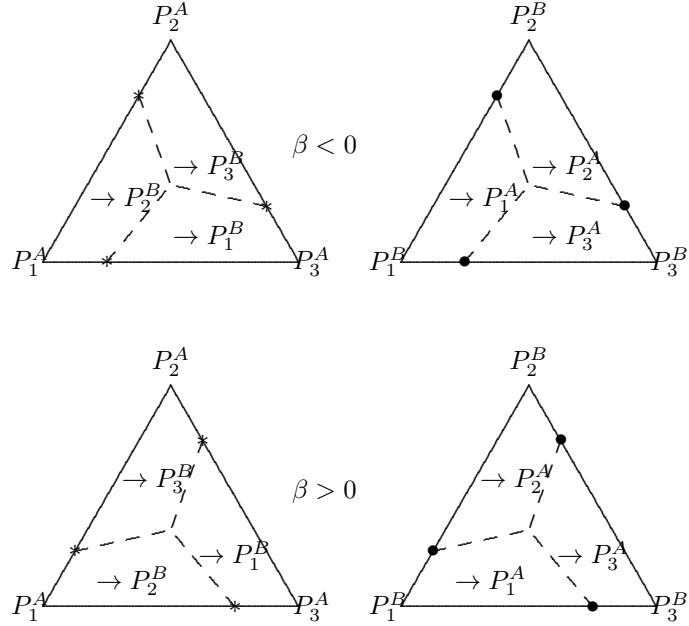
 \hfil
\beginpicture
\dimen0=0.17cm
\setcoordinatesystem units <\dimen0,\dimen0>
\setplotarea x from 0 to 20, y from 0 to 17
\setlinear
\plot 0 0 20 0 /
\plot 0 0 10 17.3 20 0 /
\setdashes
\plot 5 0 10 6 /
\plot 7.5 13 10 6 /
\plot 17.5 4.33 10 6 /
\put {$*$} at 5 0
\put {$*$} at 7.5 13
\put {$*$} at 17.5 4.33
\put {$P^A_1$} at -1 0
\put {$P^A_2$} at 10 19.2
\put {$P^A_3$} at 21 0
\put {$\to P^B_2$} at 6.5 5
\put {$\to P^B_1$} at 13 2
\put {$\to P^B_3$} at 13 7.5
\put {$\beta<0$} at 22 9
\setcoordinatesystem units <\dimen0,\dimen0> point at -28 0
\setplotarea x from 0 to 30, y from 0 to 17
\setsolid
\setlinear
\plot 0 0 20 0 /
\plot 0 0 10 17.3 20 0 /
\setdashes
\plot 5 0 10 6 /
\plot 7.5 13 10 6 /
\plot 17.5 4.33 10 6 /
\put {$\bullet$} at 5 0
\put {$\bullet$} at 7.5 13
\put {$\bullet$} at 17.5 4.33
\put {$P^B_1$} at -1 0
\put {$P^B_2$} at 10 19.2
\put {$P^B_3$} at 21 0
\put {$\to P^A_1$} at  6.5 5
\put {$\to P^A_3$} at 13 2
\put {$\to P^A_2$} at 13 7.5
\setcoordinatesystem units <\dimen0,\dimen0> point at 0 27
\setplotarea x from 0 to 30, y from -5 to 10
\put {$\beta>0$} at 22 9
\setlinear
\setsolid
\plot 0 0 20 0 /
\plot 0 0 10 17.3 20 0 /
\setdashes
\plot 15 0 10 6 /
\plot 2.5 4.33 10 6 /
\plot 12.5 13 10 6 /
\put {$*$} at 15 0
\put {$*$} at 2.5 4.33
\put {$*$} at 12.5 13
\put {$P^A_1$} at -1 0
\put {$P^A_2$} at 10 19.2
\put {$P^A_3$} at 21 0
\put {$\to P^B_2$} at 7 2
\put {$\to P^B_1$} at 15 4
\put {$\to P^B_3$} at 9 10
\setcoordinatesystem units <\dimen0,\dimen0> point at -28 27
\setplotarea x from 0 to 30, y from -5 to 10
\setsolid
\setlinear
\plot 0 0 20 0 /
\plot 0 0 10 17.3 20 0 /
\setdashes
\plot 15 0 10 6 /
\plot 2.5 4.33 10 6 /
\plot 12.5 13 10 6 /
\put {$\bullet$} at 15 0
\put {$\bullet$} at 2.5 4.33
\put {$\bullet$} at 12.5 13
\put {$P^B_1$} at -1 0
\put {$P^B_2$} at 10 19.2
\put {$P^B_3$} at 21 0
\put {$\to P^A_1$} at 7 2
\put {$\to P^A_3$} at 15 4
\put {$\to P^A_2$} at 9 10
\endpicture
\caption{\label{fig:ssb-other}
{\small The simplices $\Sigma_A$ and $\Sigma_B$ for $\beta<0$ (above) and $\beta>0$
(below). 
For $\beta\in [-1,1]$,
the sets where Player A is indifferent
between two strategies are line segments $Z^A_{ij}$ in $\Sigma_B$
connecting the equilibrium point $E^A=(1/3,1/3,1/3)$ to the points
$\frac{1}{(2-\beta)}(1,1-\beta,0)$
$\frac{1}{(2-\beta)}(1-\beta,0,1)$ and 
$\frac{1}{(2-\beta)}(0,1,1-\beta)$
(marked with $\bullet$).
Player B is indifferent between two strategies
along line segments $Z^B_{i,j}$ in $\Sigma_A$
connecting the midpoint $(1/3,1/3,1/3)$ to the points
$\frac{1}{(2+\beta)}(1+\beta,1,0)$, 
$\frac{1}{(2+\beta)}(1,0,1+\beta)$ and 
$\frac{1}{(2+\beta)}(0,1+\beta,1)$
(marked with $*$).
For $\beta=0$ these points lie on the middle
of the sides of the triangles $\partial \Sigma_A$ and
$\partial \Sigma_B$ (as in Figure \ref{fig:ssb=0}).
Increasing $\beta\in [-1,1]$ rotates the lines anticlockwise
(so the end point moves for example from $P^B_2$ towards $P^B_1$).
Note that these indifference lines in the simplices $\Sigma_A$ and $\Sigma_B$
do not coincide (despite the appearance of the figure); 
for example, when $\beta=1$ then in $\Sigma_B$ the
indifference lines $Z^A_{i,j}$ go through the corners,
but the lines $Z^B_{i,j}$ in $\Sigma_A$ do not.
}} 
\end{figure}

Figure~\ref{fig:ssb-other} shows pictures of the phase space
in the cases $-1 < \beta < 0$ and
$0 <\beta < 1$, marking
the lines where the players are indifferent by dashed lines. 
(In Figure~\ref{fig:ssb=0} the phase space was drawn for $\beta=0$.)
Note that for all values of $\beta$ both players are 
indifferent between all three strategies
at the point $E=(E^A, E^B)$ where  $E^A=(E^B)^T=(1/3,1/3,1/3)$.
Note that when $\beta\ne 0$,
for each $p^A\in  Z_{ij}^B\subset \Sigma_A$, 
each aim vector $p^AP_k^A$ is transversal to $Z_{ij}^B$,
and similarly for $p^B\in Z_{ij}^A \subset \Sigma^B$.
So the transversality condition of \ref{eqn:trans} holds
for $\beta\ne 0$.

\subsection{Continuity and uniqueness of the flow}
\label{ss:ABcont}

We now apply the results of the previous section to the
family of examples defined by matrices (\ref{eqn:AB}).
An essential difference between
the cases $\beta<0$ and $\beta>0$ can be understood in terms 
of that analysis. 

\begin{prop}[Continuity and uniqueness] \label{initre}
When $\beta \in (0,1)$ the unique flow on the complement
of the indifference set $Z$ extends to a flow on $Z$ 
which is everywhere continuous and
unique except at $E =(E^A,E^B)$.
When $\beta\in(-1,0]$, the unique flow on the complement
of $Z$ does not extend continuously to some codimension
two surfaces in $Z^*$. In more detail:

\begin{itemize}
\item If $\beta\in (-1,0)$ then 
the transversality condition
from Proposition \ref{contitra} 
holds for all
codimension one indifference points in $Z$.
All the codimension two points in $Z^*$ are in cases (1) or (3), and
the flow off $Z$ extends continuously to those points that are in case (1).
The differential inclusion (\ref{eq1}) defines a flow that 
is unique except at those points, where the flow is discontinuous.
\item If $\beta=0$ then the transversality condition
from Proposition \ref{contitra} 
fails to hold
in several codimension-one sets (for example,
the set $Z^B_{1,2}\times \{z\in \Sigma_B\st
BR_A(p^B)=P^A_1\}$). 
\item If $\beta\in (0,1)$ then the transversality condition
from Proposition \ref{contitra} 
holds on
all codimension one indifference sets,
and all the codimension two points in $Z^*$
are in case (1) or (2). In particular, the flow off $Z$ extends continuously 
to all points except $E=(E^A, E^B)$, and so the differential
inclusion (\ref{eq1}) defines a unique and continuous flow everywhere
except at the interior equilibrium $E$.
\end{itemize}
\end{prop}
\begin{proof}
This can be seen just by
looking at the possible $2 \times 2$ games that
can be obtained from submatrices of $A$ and $B$ and applying the
results from the previous section. (It is also possible to
convince oneself of the result graphically by 
comparing the disposition of lines $Z^{A/B}_{i,j}$
and points $P^A_i$ and $P^B_j$ in Figure ~\ref{fig:ssb=0} and
\ref{fig:ssb-other}.)
When $\beta<0$ case (3) holds, for example, in $Z^B_{12} \times
Z^A_{2,3}$.
When $\beta>0$
there are points in $Z^*$, for example on the set
$Z^B_{1,2} \times Z^A_{1,2}$ where case (2)
arises and trajectories locally spiral around and
tend towards $Z^*$ (until they meet another set in $Z^*$).
On the other hand orbits cross $Z^B_{1,2}\times Z^A_{2,3}$
(which is in case (1) transversally when $\beta>0$. 
\end{proof}

We will see later that when $\beta>0$ there is genuine non-uniqueness
of the flow at $E$.
The spiralling behaviour that can occur near $Z^*$ for $\beta>0$
contributes significantly to the extra richness of behaviour in $\beta>0$.
The subset 
$J \subset Z^*$ 
\footnote{We could take $J$ to stand for "jitter", 
since as nearby trajectories
spiral around $J$, both players jitter back and forth between strategies.} on which we have this behaviour
consists of  the six pieces of $Z^*$:
$$J=
(Z^B_{1,2}\times Z^A_{3,1})\cup (Z^B_{1,2}\times Z^A_{1,2})
\bigcup
(Z^B_{2,3}\times Z^A_{1,2})\cup (Z^B_{1,2}\times Z^A_{2,3})\bigcup
$$
\begin{equation}
\bigcup
(Z^B_{3,1}\times Z^A_{2,3})\cup (Z^B_{2,3}\times Z^A_{3,1})
\label{eqn:defJ}
\end{equation}
The set consisting of the remaining three pieces of $Z^*$:
\begin{equation}
 T= (Z^B_{1,2}\times Z^A_{2,3}) \bigcup (Z^B_{2,3}\times Z^A_{3,1})
\bigcup (Z^B_{3,1}\times Z^A_{1,2})
\label{eqn:defT}
\end{equation}
contains the codimension two points where the flow is transversal (case (1)).

\subsection{Existence and stability of a clockwise Shapley periodic orbit for
$\beta \in (-1,\sigma)$}

Let $\sigma$ be the golden mean.
We shall first show the existence of a stable periodic orbit 
for $\beta \in (-1, \sigma)$.  We call this the Shapley orbit, as it is the continuation
of the  orbit described by \citeasnoun{Shapley64} for the case $\beta=0$, 
or the clockwise orbit, 
as the projections of the orbit into the spaces $\Sigma_A$ and $\Sigma_B$
rotate clockwise about the equilibrium, see the figure below.
In the next subsection we shall also show that this orbit is globally attracting
when $\beta\le 0$.
For $\beta=0$ orbits spiral away from the equilibrium points
towards a periodic orbit (which consists in 4 dimensional 
phase space $\Sigma_A\times \Sigma_b$ of 6 pieces 
of straight lines).


\includegraphics{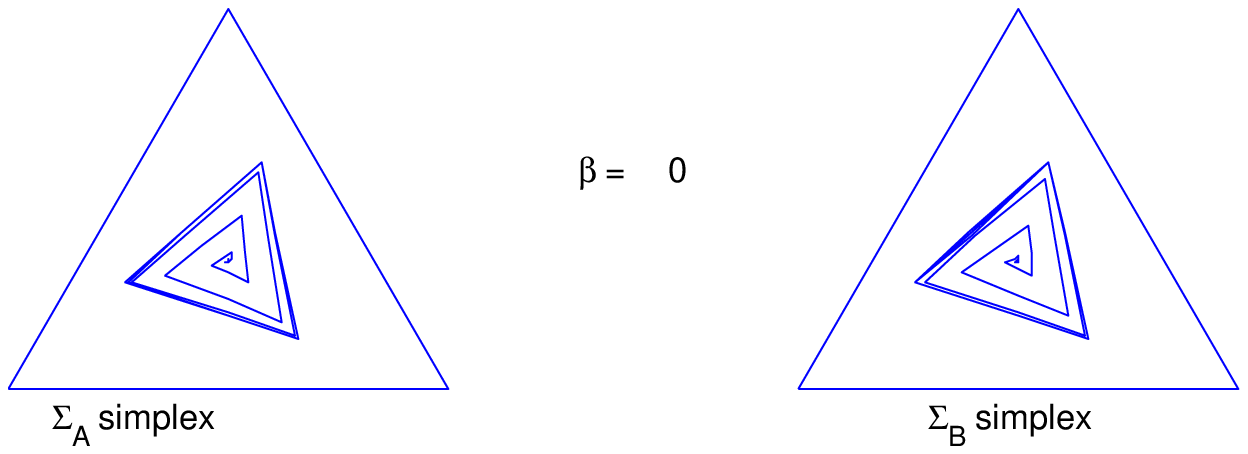}

The figure below shows the utility $n=Ap^B$ of player $A$
as a function of time for $\beta=0$. The players move away
(in an oscillatory fashion)  from the equilibrium point
and tend to a periodic motion. Note that where two graphs
meet, the players are indifferent between two strategies.

\includegraphics[width=4in]{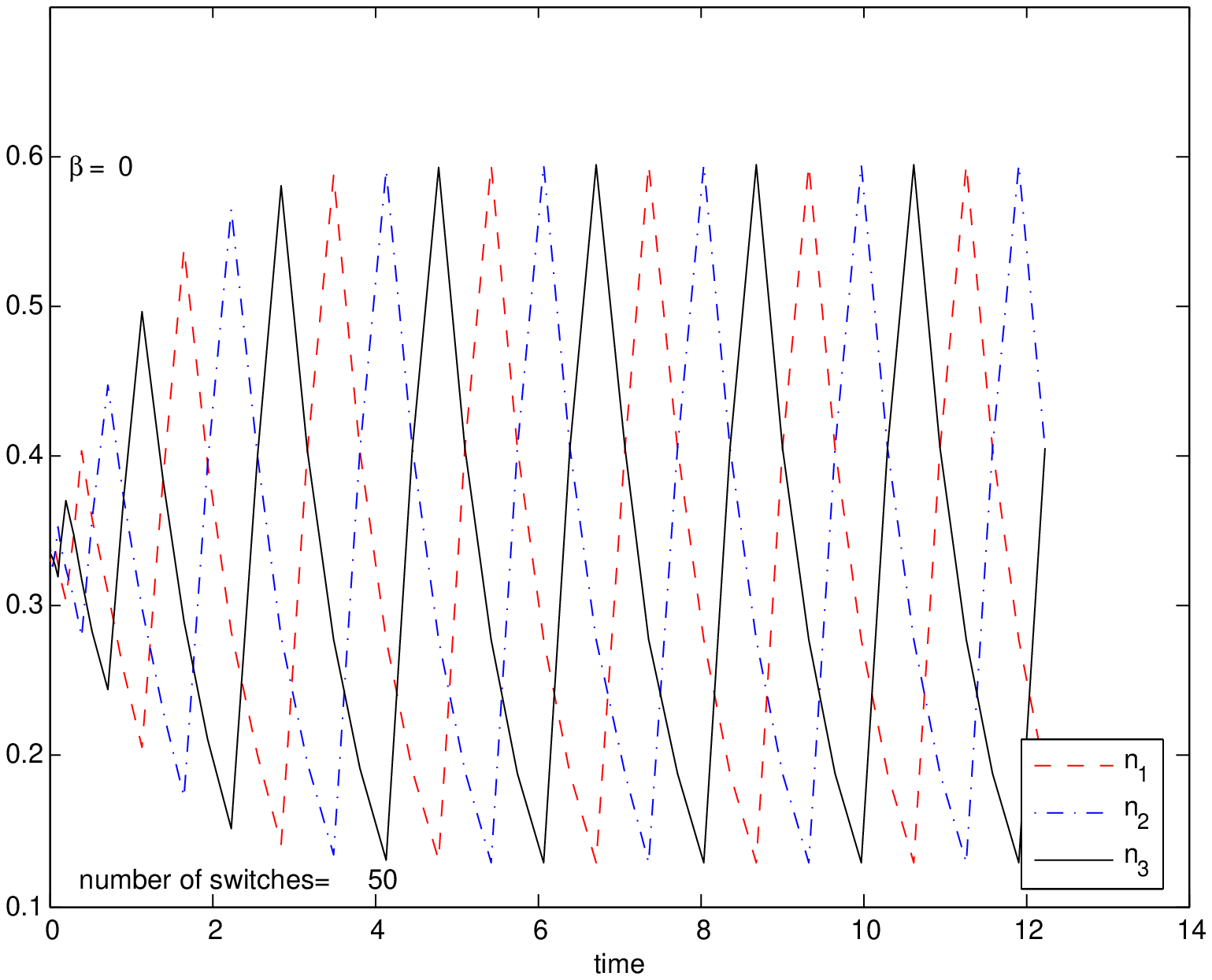}

\begin{theo}[Existence of the Shapley periodic orbit]
\label{exipeshap}
For $\beta\in (-1,\sigma)$ there exists a periodic orbit in which
player 1 and 2 infinitely repeat the strategies
$$ 
\begin{array}{cccccccc}	
\mbox{ time \quad } & 1 & 2 & 3 & 4 & 5 & 6  &  \\
\mbox{strategy player A \quad }  &1 & 2 & 2 & 3  & 3 & 1 & \\ 
\mbox{strategy player B \quad }  & 2 & 2 & 3 & 3 & 1 & 1 & \end{array}
$$
(So player A copies B and player B tries to play one ahead of $A$.)
This is a continuous deformation of the  orbit that Shapley found for $\beta=0$, 
and as $\beta\uparrow \sigma$
this orbit shrinks in diameter to zero and tends to the point
$(E^A, E^B)$. For $\beta\in (\sigma,1)$ this periodic 
orbit no longer exists.
\end{theo}

The proof of this Theorem is given in Section~\ref{ss:appclock}
(in the Appendix).
Figure~\ref{figbis0} (and later figures up to Figure~\ref{figbispos}) 
allow us to illustrate schematically the relationship between 
regions $S_{ij}$ where player A prefers strategy $i$ and Player
B prefers strategy $j$, $1 \le i,j \le 3$. The 9 squares in these 
diagrams each represent one such region, as labelled. We think of 
the $3 \times 3$ checker board as a torus, so that the top and bottom are
identified, as are the left and right sides (so if we go for a walk 
on the checkerboard and wander off the top, we return with the same 
horizontal coordinate at the bottom, and similarly if we wander off
the left of the diagram we return with the same vertical coordinate 
on the right). With this convention, the diagram illustrates the  
adjacency of regions; if just one
player changes her preferred strategy, we cross a horizontal or 
vertical line (representing a codimension-one set) to get to the 
new region. If we wish both players to change their preferred strategies 
at the same instant, we must cross diagonally through a corner 
(representing a codimension-two set) into a new region. 
Higher codimension points, such as the internal equilibrium $(E^A, E^B)$
where all regions meet, are not represented. Note also that
each region is actually 4-dimensional, rather than 2 as represented 
in the diagram, and so if we wish to project trajectories of the system 
onto the diagram we may find that they cross; certainly we do not want to
claim there is any sort of unique flow determined on the checkerboard.  

The usefulness 
of these diagrams is that we may add arrows showing which 
transitions are allowed, which we will use in the proof of the theorem. 
Thus, for example, if we imagine that we begin in the top left-hand 
region where $A \to P^A_1$ and $B \to P^B_2$
then inspection of Figure~\ref{fig:ssb=0} for $\beta=0$ shows that 
the only region we can move to (and indeed will move to in due course) 
is the region where $A \to P^A_2$ and $B \to P^B_2$, 
as represented by the arrow 
between these two regions in Figure \ref{figbis0}.
The other arrows are computed similarly. In regions such 
as $A \to P^A_1$ and $B \to P^B_3$ there is a codimension-one set such that 
trajectories started on this set will reach indifference lines
in $\Sigma_A$ and $\Sigma_B$ simultaneously. This set is represented 
by the dotted line in the appropriate square in Figure \ref{figbis0}, 
and similarly in other squares where this occurs. Finally, some
of the codimension one indifference surfaces are represented by dashed 
lines in the $\beta=0$ case.  This
is rather special, and is where the 
transversality condition of Proposition \ref{contitra} fails; orbits 
move off in different directions on opposite sides of these lines. 
Figure \ref{figbisneg} for $\beta<0$ is
rather similar, and the extra dotted lines are explained in the caption.

\begin{figure}[htp]
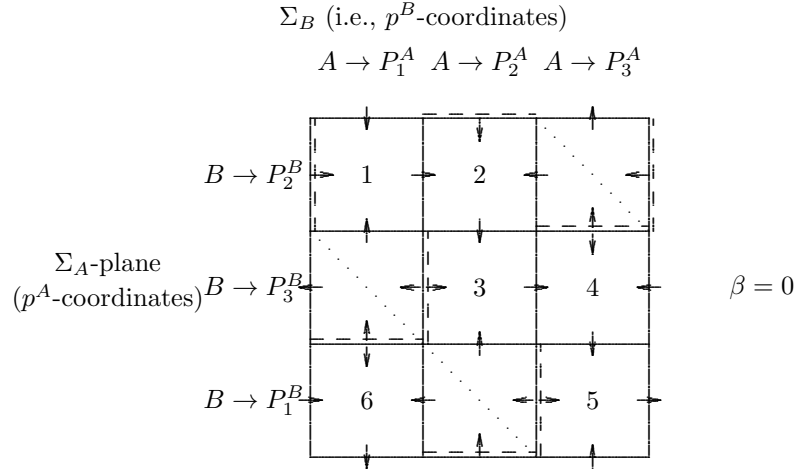
 \hfil
\beginpicture
\dimen0=0.15cm
\setcoordinatesystem units <\dimen0,\dimen0>
\setplotarea x from 0 to 30, y from 0 to 42
\setlinear
\put {$\Sigma_B$ (i.e., $p^B$-coordinates)} at 10 39
\put {$A \to P^A_1$} at 5 35
\put {$A \to P^A_2$} at 15 35
\put {$A \to P^A_3$} at 25 35
\put {$B \to P^B_2$} at -5 25
\put {$B \to P^B_3$} at -5 15  \put {$\beta=0$} at 40 15
\put {$B \to P^B_1$} at -5 5   \put {$\Sigma_A$-plane} at -18 17
                           \put {($p^A$-coordinates)} at -18 14
\arrow <5pt> [0.2,0.4] from 0 25 to 2 25
\arrow <5pt> [0.2,0.4] from 30 25 to 28 25
\arrow <5pt> [0.2,0.4] from 25 20.3 to 25 22
\arrow <5pt> [0.2,0.4] from 25 19.7 to 25 18
\arrow <5pt> [0.2,0.4] from 15 30 to 15 28
\arrow <5pt> [0.2,0.4] from 5 10.3 to 5 12
\arrow <5pt> [0.2,0.4] from 5 9.7 to 5 8
\arrow <5pt> [0.2,0.4] from 9.7 15 to 8 15
\arrow <5pt> [0.2,0.4] from 10.3 15 to 12 15
\arrow <5pt> [0.2,0.4] from 15 0 to 15 2
\arrow <5pt> [0.2,0.4] from 19.7 5 to 18 5
\arrow <5pt> [0.2,0.4] from 20.3 5 to 22 5
\put {$1$} at 5 25
\put {$2$} at 15 25
\put {$3$} at 15 15
\put {$4$} at 25 15
\put {$5$} at 25 5
\put {$6$} at 5 5
\plot 0 0 30 0 /
\plot 0 10 30 10 /
\plot 0 20 30 20 /
\plot 0 30 30 30 /
\plot 0 0 0 30  /
\plot 10 0 10 30  /
\plot 20 0 20 30  /
\plot 30 0 30 30  /
\setdots
\plot 0 20 10 10 /
\plot 10 10 20 0 /
\plot 20 30 30 20 /
\setsolid
\arrow <5pt> [0.2,0.4] from 5 1 to 5 -1
\setdashes \plot 0 10.4 10 10.4 / \setsolid
\arrow <5pt> [0.2,0.4] from 5 19 to 5 21
\arrow <5pt> [0.2,0.4] from 5 31 to 5 29
\arrow <5pt> [0.2,0.4] from 15 9 to 15 11
\arrow <5pt> [0.2,0.4] from 15 21 to 15 19
\setdashes \plot 10 0.4 20 0.4 / \setsolid
\setdashes \plot 10 30.4 20 30.4 /\setsolid
\arrow <5pt> [0.2,0.4] from 25 -1 to 25 1
\arrow <5pt> [0.2,0.4] from 25 11 to 25 9
\setdashes \plot 20 20.4 30 20.4 /\setsolid
\arrow <5pt> [0.2,0.4] from 25 29 to 25 31
\arrow <5pt> [0.2,0.4] from 29 5 to 31 5
\arrow <5pt> [0.2,0.4] from 11 5 to 9 5
\arrow <5pt> [0.2,0.4] from -1 5 to 1 5
\setdashes \plot 20.4 0 20.4 10 /\setsolid
\arrow <5pt> [0.2,0.4] from 31 15 to 29 15
\arrow <5pt> [0.2,0.4] from 19 15 to 21 15
\arrow <5pt> [0.2,0.4] from 1 15 to -1 15
\setdashes \plot 10.4 10 10.4 20 /\setsolid
\arrow <5pt> [0.2,0.4] from 21 25 to 19 25
\arrow <5pt> [0.2,0.4] from 9 25 to 11 25
\setdashes \plot 0.4 20 0.4 30 /\setsolid
\setdashes \plot 30.4 20 30.4 30 /\setsolid
\endpicture
\caption{\label{figbis0}
{\small A diagram showing allowed transitions between regions of pure preference
for $\beta=0$. See the text for explanations of dotted and dashed lines.}
}
\end{figure}

\begin{figure}[htp]
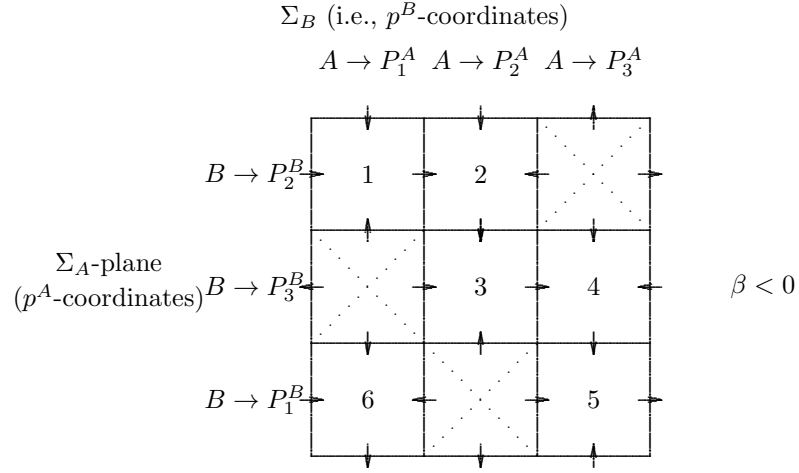
 \hfil
\beginpicture
\dimen0=0.15cm
\setcoordinatesystem units <\dimen0,\dimen0>
\setplotarea x from 0 to 30, y from 0 to 30
\setlinear
\put {$\Sigma_B$ (i.e., $p^B$-coordinates)} at 10 39
\put {$A \to P^A_1$} at 5 35
\put {$A \to P^A_2$} at 15 35
\put {$A \to P^A_3$} at 25 35
\put {$B \to P^B_2$} at -5 25
\put {$B \to P^B_3$} at -5 15  \put {$\beta<0$} at 40 15
\put {$\Sigma_A$-plane} at -18 17
                           \put {($p^A$-coordinates)} at -18 14
\put {$B \to P^B_1$} at -5 5   
\put {$1$} at 5 25
\put {$2$} at 15 25
\put {$3$} at 15 15
\put {$4$} at 25 15
\put {$5$} at 25 5
\put {$6$} at 5 5
\plot 0 0 30 0 /
\plot 0 10 30 10 /
\plot 0 20 30 20 /
\plot 0 30 30 30 /
\plot 0 0 0 30  /
\plot 10 0 10 30  /
\plot 20 0 20 30  /
\plot 30 0 30 30  /
\setdots
\plot 0 20 10 10 / \plot 0 10 10 20 /
\plot 10 10 20 0 / \plot 10 0 20 10 /
\plot 20 30 30 20 / \plot 20 20 30 30 /
\setsolid
\arrow <5pt> [0.2,0.4] from 5 1 to 5 -1
\arrow <5pt> [0.2,0.4] from 5 11 to 5 9
\arrow <5pt> [0.2,0.4] from 5 19 to 5 21
\arrow <5pt> [0.2,0.4] from 5 31 to 5 29
\arrow <5pt> [0.2,0.4] from 15 9 to 15 11
\arrow <5pt> [0.2,0.4] from 15 21 to 15 19
\arrow <5pt> [0.2,0.4] from 15 1 to 15 -1
\arrow <5pt> [0.2,0.4] from 15 31 to 15 29
\arrow <5pt> [0.2,0.4] from 25 -1 to 25 1
\arrow <5pt> [0.2,0.4] from 25 11 to 25 9
\arrow <5pt> [0.2,0.4] from 25 21 to 25 19
\arrow <5pt> [0.2,0.4] from 15 21 to 15 19
\arrow <5pt> [0.2,0.4] from 25 29 to 25 31
\arrow <5pt> [0.2,0.4] from 29 5 to 31 5
\arrow <5pt> [0.2,0.4] from 11 5 to 9 5
\arrow <5pt> [0.2,0.4] from -1 5 to 1 5
\arrow <5pt> [0.2,0.4] from 19 5 to 21 5
\arrow <5pt> [0.2,0.4] from 31 15 to 29 15
\arrow <5pt> [0.2,0.4] from 19 15 to 21 15
\arrow <5pt> [0.2,0.4] from 1 15 to -1 15
\arrow <5pt> [0.2,0.4] from 9 15 to 11 15
\arrow <5pt> [0.2,0.4] from 21 25 to 19 25
\arrow <5pt> [0.2,0.4] from -1 25 to 1 25
\arrow <5pt> [0.2,0.4] from 9 25 to 11 25
\arrow <5pt> [0.2,0.4] from 29 25 to 31 25
\endpicture
\caption{\label{figbisneg}
{\small The transition diagram of the
flow for $\beta<0$. The dotted lines again
indicate the codimension-one sets
where $p^A$ and $p^B$ leave the relevant regions
simultaneously.
The points where these dotted lines meet,
correspond to a set of initial conditions where the trajectory exits
the region 
by hitting $E^A=E^B=(1/3,1/3,1/3)$.} }
\end{figure}

Let us call the discontinuity set the subset of $Z^*$ for which
type (1) behaviour occurs. It has (of course) measure zero.

\begin{theo} [The Shapley orbit is globally attracting for $\beta\le 0$]
For $\beta\le 0$, the flow has a unique closed orbit (the Shapley orbit);
the flow of all initial conditions which do not enter
the discontinuity set tends to this closed orbit. \label{thm:shapley}
\end{theo}
\begin{proof}
When $\beta\le 0$, then the only regions which can be entered are
those marked with a central integer label. So any trajectory which does not
enter the discontinuity set, eventually enters one of these
regions. (These orbits form a set of zero Lebesgue measure.)
From the diagram it is clear that from then on
the orbit cyclically visits the regions 1-6.
Now consider, for example, the part $S_1$ of the space $Z^A_{1,2}$ 
where player A is indifferent between strategies $P^A_2$ and $P^A_1$, 
and player B strictly prefers strategy $P^B_2$ (this corresponds in the diagram 
to the line dividing regions labelled 1 and 2).
Similarly consider the region $S_2$ where the orbits leave
region labelled 2.
Each of these sets $S_i$ is the product of a line segment
and a kite-shaped quadrangle, with one boundary point 
the equilibrium point $E$. 
Now consider the map $T$ which assigns to a point in
 $S_1$ the first visit of the forward orbit through $x$ to  $S_2$. Since all orbits in region 1 
will aim for the same point and only exit via $S_2$ (when $\beta\in (-1,0]$), this map is continuous;
moreover, it has the property that
$T(p^A,p^B)\in \mbox{interior}(S_2)$
for each  $(p^A,p^B)\in \overline{S_1}\setminus  
(\overline{S_1}\cap \overline{S_2})$. In fact, $T$ is a central projection.
In other words, if we put linear coordinates on $S_1$ and $S_2$
then $T\colon S_1\to S_2$ takes the form
\begin{equation}
T(z)=\frac{1}{f(z)}L(z)
\end{equation}
where $L\colon S_1\to S_2$ and $f\colon S_1 \to \rz$
are affine transformations.
So in linear coordinates,
$$T(x,y,z)=\left(
\frac{a_1+b_1x+c_1y+d_1z}
{a+b x + c y + d z},
\frac{a_2+b_2x+c_2y+d_2z}
{a+b x + c y + d z},
\frac{a_3+b_3x+c_3y+d_3z}
{a+b x + c y + d z}\right).$$
(These maps are explicitly calculated in the appendix; here
we shall only use that it is of this form.)
Note that compositions of such maps are again of this form.
If
$$T_1=\frac{1}{f_1}L_1\mbox{ and }T_2=\frac{1}{f_2}L_2$$
are two such transformations and
we take linear coordinates such that $T_1(0)=0$ (i.e.
$L_1(0)=0$)
then the composition formula becomes particularly simple:
\begin{equation}
T_2\circ T_1=\frac{1}{f_2(T_1(z))\cdot f_1(z)} L_2\circ L_1(z).
\label{compositionform}
\end{equation}
Note that because of cancellations 
$f_2(T_1(z))\cdot f_1(z)$ actually is an affine map in $z$.

Now let $P$ be the first return map to $S_1$. 
As we observed $P$ is a continuous central projection as above
(it is composition of such maps),
so of the form $z\mapsto P(z)=\frac{1}{f(z)}L(z)$ with
$P(S_1)\subset S_1$ and, because 
$T\colon S_1\to S_2$ maps  all points except those in
$\overline{S_1}\cap \overline{S_2}$ into $\mbox{interior}(S_2)$ 
(and similarly for the other transition maps),
\begin{equation}
P(\overline{S_1})\subset \mbox{interior}(S_1)\cup E .
\label{eq:almoststrict}
\end{equation}
In addition, $P$ has a
fixed $p$ in the interior of $S_1$ corresponding to the Shapley periodic orbit.
Let us show that $\lim_{n\to \infty}P^n(x)\to p$ for each $x\in S_1\setminus E$,
where $P^n$ is the $n$-th iterate of $P$.
To see this, notice that  $P$ sends 
lines to lines
and restricted to each line, 
$T$ is a Moebius transformation, i.e. of the form 
$t\mapsto \frac{\alpha+t \beta}{\gamma+t\delta}$
(where we take affine coordinates on the lines).
Let us first consider the line $l_0$ through $p$ and $E$,
and let $b$ be so that $l_0\cap S_1=[E,b]$.
This line $l_0$ is mapped into itself because $p,E\in l_0$
are both fixed points.  
$P\colon l_0\to l_0$  is a Moebius
transformation which sends the line segment
$[E,b]=l_0\cap S_1$ continuously  into $[a,b)$
(this is because of (\ref{eq:almoststrict}))
and having fixed point $p\in (E,b)$.
It follows from this that for any $x\in (E,b]$,
$\lim_{n\to \infty}P^n(x)\to p$.
Now consider any other eigendirection $l$ through $x$.
Because of (\ref{eq:almoststrict}), $P$ maps 
$l\cap S_1$ continuously strictly inside itself, and since 
the restriction of  $P$ to this line is a Moebius transformation,
again all points in $l\cap S_1$ converge to $p$.
It $P$ has a complex eigenvalue then the argument is similar.

To consider any other point $x$, 
take the projective space $\PP$
of lines through $p$ in  $S_1$.  The return map $P$
assigns to each line $\hat x\in \PP$ a new line $\hat P(\hat x)\in \PP$.
The eigenvalues of the linear transformation 
$L$ correspond to fixed points of $\hat P$,
and for each line $l$ through $x$ in $S_1$, 
$P^n(l)$ tends to an eigenspace of the
linearisation of $P$ at $x$. Combining this with the
above, it follows that for each $x\in S_1\setminus E$,
$P^n(x)\to p$.

It follows that for the original flow
all trajectories which do not enter the discontinuity set
are attracted to the periodic Shapley orbit.
\end{proof}

This theorem completes our investigation of $\beta \le 0$.

\subsection{Existence and Stability of a Shapley and also
an anti-Shapley periodic orbit for $\beta > 0$}

\begin{figure}[htp]
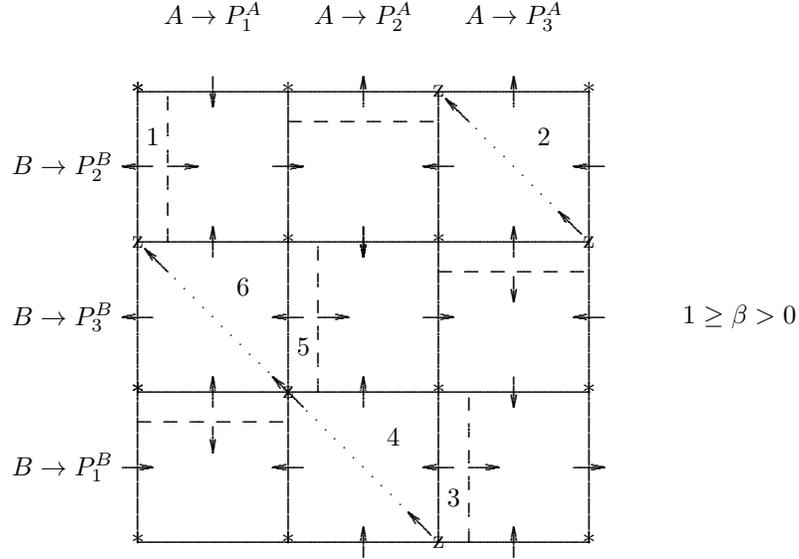
 \hfil
\beginpicture
\dimen0=0.2cm
\setcoordinatesystem units <\dimen0,\dimen0>
\setplotarea x from 0 to 30, y from 0 to 30
\setlinear
\put {$A \to P^A_1$} at 5 35
\put {$A \to P^A_2$} at 15 35
\put {$A \to P^A_3$} at 25 35
\put {$B \to P^B_2$} at -5 25
\put {$B \to P^B_3$} at -5 15  \put {$1\ge \beta>0$} at 40 15
\put {$B \to P^B_1$} at -5 5
\put {$1$} at 1 27
\put {$2$} at 27 27
\put {$3$} at 21 3
\put {$4$} at 17 7
\put {$5$} at 11 13
\put {$6$} at  7 17
\put {*} at 0 30
\put {*} at 30 0
\put {*} at 10 0
\put {*} at 10 30
\put {*} at 30 30
\put {*} at 10 20
\put {*} at 20 20
\put {*} at 0 10
\put {*} at 20 10
\put {*} at 30 10
\put {*} at 0 0
\put {*} at 0 10
\put {*} at 0 30
\put {z} at 20 30
\put {z} at 0 20
\put {z} at 30 20
\put {z} at 10 10
\put {z} at 20 0
\arrow <5pt> [0.2,0.4] from 22 5 to 24 5
\arrow <5pt> [0.2,0.4] from 12 15 to 14 15
\arrow <5pt> [0.2,0.4] from 2 25 to 4 25
\plot 0 0 30 0 /
\plot 0 10 30 10 /
\plot 0 20 30 20 /
\plot 0 30 30 30 /
\plot 0 0 0 30  /
\plot 10 0 10 30  /
\plot 20 0 20 30  /
\plot 30 0 30 30  /
\setdots
\plot 0 20 10 10 /
\plot 10 10 20 0 /
\plot 20 30 30 20 /
\setsolid
\arrow <5pt> [0.2,0.4] from 22 28 to 20.5 29.5
\arrow <5pt> [0.2,0.4] from 29.5 20.5 to 28 22
\arrow <5pt> [0.2,0.4] from 19.5 0.5 to 18 2
\arrow <5pt> [0.2,0.4] from 11 9 to 9 11
\arrow <5pt> [0.2,0.4] from 2 18 to 0.5 19.5
\setdashes \plot 0 8 10 8 / \setsolid
\arrow <5pt> [0.2,0.4] from 5 9 to 5 11
\arrow <5pt> [0.2,0.4] from 9 25 to 11 25
\arrow <5pt> [0.2,0.4] from 5 7.7 to 5 6
\arrow <5pt> [0.2,0.4] from 5 19 to 5 21
\arrow <5pt> [0.2,0.4] from 5 31 to 5 29
\arrow <5pt> [0.2,0.4] from 15 21 to 15 19
\setdashes \plot 10 28 20 28 / \setsolid
\arrow <5pt> [0.2,0.4] from 15 -1 to 15 1
\arrow <5pt> [0.2,0.4] from 15 29 to 15 31
\arrow <5pt> [0.2,0.4] from 25 -1 to 25 1
\arrow <5pt> [0.2,0.4] from 25 11 to 25 9
\setdashes \plot 20 18 30 18 /\setsolid
\arrow <5pt> [0.2,0.4] from 25 17.7 to 25 16
\arrow <5pt> [0.2,0.4] from 25 19 to 25 21
\arrow <5pt> [0.2,0.4] from 15 21 to 15 19
\arrow <5pt> [0.2,0.4] from 25 29 to 25 31
\arrow <5pt> [0.2,0.4] from 29 5 to 31 5
\arrow <5pt> [0.2,0.4] from 11 5 to 9 5
\arrow <5pt> [0.2,0.4] from 15 9 to 15 11
\arrow <5pt> [0.2,0.4] from -1 5 to 1 5
\setdashes \plot 22 0 22 10 /\setsolid
\arrow <5pt> [0.2,0.4] from 21 5 to 19 5
\arrow <5pt> [0.2,0.4] from 31 15 to 29 15
\arrow <5pt> [0.2,0.4] from 19 15 to 21 15
\arrow <5pt> [0.2,0.4] from 1 15 to -1 15
\setdashes \plot 12 10 12 20 /\setsolid
\arrow <5pt> [0.2,0.4] from 11 15 to 9 15
\arrow <5pt> [0.2,0.4] from 21 25 to 19 25
\setdashes \plot 2 20 2 30 /\setsolid
\arrow <5pt> [0.2,0.4] from 1 25 to -1 25
\arrow <5pt> [0.2,0.4] from 31 25 to 29 25
\endpicture
\caption{\label{figbispos}
{\small The transition diagram 
of the flow for $\beta\in (0,1)$. Boxes again represent regions
$S_{ij}$, and arrows indicate possible transitions between
regions. Routes may now spiral around
near the `corners' marked with $*$; these 6 corners represent 
codimension-two points on the six pieces of $J$. The  3 corners 
marked $z$ represent codimension-two points on the 3 pieces 
of $Z^* \setminus J$ and trajectories cross these sets transversally;
there is no spiralling. 
The numbers $1-6$ in the regions
show a possible route for an anticlockwise periodic orbit. Notice that the
route for the clockwise Shapley orbit is still also allowed.
The dashed lines correspond to preimages
of points $(E^A,p^B)$  or $(p^A, E^B)$,
whereas the dotted lines correspond to preimages
of $Z^* \setminus J$, but where the orbits will cross transversally 
as in case (1) above. When $\beta=1$
the situation remains almost the same except that it is no longer possible
to have two types
of exits from the regions with the dotted lines (no exit is possible in the
horizontal direction in the three relevant boxes, but this restriction does
not prevent the continued existence of the anticlockwise route marked).}}
\end{figure}

Let us now consider orbits when $\beta\in (0,1)$.
The transition diagram for this case is given in
Figure~\ref{figbispos} and it is immediately clear that there are now many
more possible routes that may be taken from region to region. We stress that
the existence of a possible route allowed by the arrows, does not
necessarily imply the existence of a corresponding trajectory in the flow; 
we will need to perform some calculations to see which routes are actually 
realised for any given value of $\beta$. In this figure the dashed lines 
represent codimension-one sets of points from which trajectories tend, 
for one or other player, to $E^A$ or $E^B$. The sequence of regions
visited will be different depending which side of these sets you begin.
[We should note that the dashed curves
do not move entirely across the squares as $\beta \to 1$ (because
the points in Figure~\ref{fig:ssb-other} marked with $\bullet$ in $\Sigma_B$
still have not reached the corners of $\Sigma_B$ when $\beta=1$).]

We show in the Appendix, by direct calculation, 
that precisely one of the anticlockwise 
or clockwise orbits mentioned in the figure caption exists 
when $\beta\ne \sigma$, where 
$\sigma=(-1+\sqrt{5})/2$,  the golden mean. In particular,
for $\beta\in [-1,\sigma)$ there is an attracting orbit which
goes round
the midpoint of $\Sigma_A$ and $\Sigma_B$ in an clockwise fashion
(this is a continuous deformation of Shapley's periodic orbit
for $\beta=0$, which we have already studied in the previous section for
$\beta \le 0$),
whereas for $\beta\in (\sigma,1)$ there is a different 
periodic orbit which goes around these points in a anticlockwise 
fashion. It is attracting for some interval of $\beta$ values to 
be described further below, and is otherwise of saddle-type.

\begin{theo}[Existence and stability of the Shapley and
anti-Shapley periodic orbits]
\label{exipe}
For $\beta\in [0,\sigma)$ there exists the periodic orbit mentioned
already in Theorem~\ref{exipeshap}, i.e. it visits,
in the same order,
the same regions as the Shapley orbit of Theorem \ref{thm:shapley}.
This is a continuous deformation of that orbit, and
as $\beta\uparrow \sigma$
this orbit shrinks in diameter to zero and tends to the point
$(E^A, E^B)$.
This orbit is stable for all $\beta\in [0,\sigma)$.

For $\beta\in (\sigma,1)$ there exists a periodic orbit (anticlockwise) which
visits the regions marked 1-6 in Figure \ref{figbispos} cyclically
in order, i.e., follow the following pattern:
$$ 
\begin{array}{cccccccc}	
\mbox{ time \quad } & 1 & 2 & 3 & 4 & 5 & 6  &  \\
\mbox{strategy player A \quad }  & 1 & 1 & 3 & 3  & 2 & 2 & \\ 
\mbox{strategy player B \quad }  & 3 & 2 & 2 & 1 & 1 & 3 & 
\end{array}
$$
As $\beta\downarrow \sigma$
this orbit shrinks in diameter to zero and tends to the point
$(E^A,E^B)$. Such a periodic orbit does NOT exist
when $\beta<\sigma$.

Shapley's (clockwise) periodic orbit is attracting for all
$\beta\in [0,\sigma)$. The anti-Shapley (anticlockwise) periodic
orbit which exists for $\beta\in (\sigma,1)$ is attracting when
$\beta\in (\tau,1)$ and of saddle-type when
$\beta\in (\sigma,\tau)$. Here $\tau$ is the root of
some specific expression and is approximately
$0.915$.
\end{theo}

We remark that it does not seem to be possible to prove the
existence of the anticlockwise orbit, even when $\beta=1$, 
by the same sort of technique as we used to prove the existence 
of the clockwise orbit in $\beta \le 0$. This is because  
even when $\beta=1$ the relevant version of Figure~\ref{figbispos}
still allows many routes through the diagram; it was only because there was
only a single route possible in the $\beta\le 0$ case that we were able 
to deduce that an orbit following that route must exist.
Neither the Shapley nor the anti-Shapley orbits
attract every starting position: as we show in the sequel
to this paper \citeasnoun{CS2006}, there are an uncountable number
of orbits which tend to the interior equilibrium point. The existence
of one such orbit is shown in the next subsection.

The proof of this Theorem is in the Appendix, as is a discussion
of the bifurcation at $\beta=\tau$. 


\subsection{Dynamics and a periodic orbit $\Gamma$ on the set $J$}

We will now consider the flow on the set $J$ (more precisely, the unique
continuous extension of the flow to $J$).
The dynamics near $J$ is very interesting
because orbits can spiral for very long time-lengths
near $J$, before moving away.
Because on  $J$ at least one of the  two players
is indifferent between two strategies, 
nearby orbits correspond to behaviour
where players frequently switch their preferred behaviour.

\begin{figure}[htp]
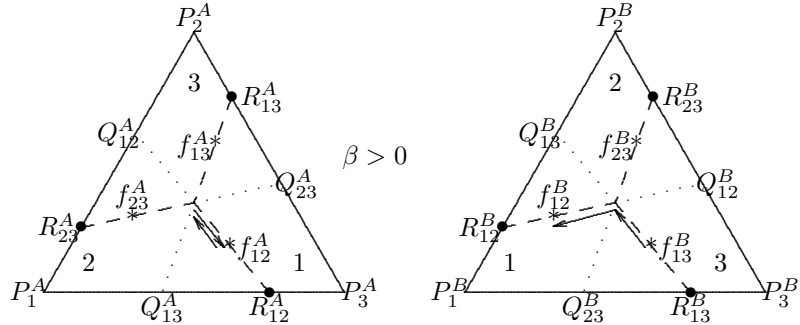
 \hfil
\beginpicture
\dimen0=0.2cm
\setcoordinatesystem units <\dimen0,\dimen0>
\setplotarea x from 0 to 30, y from -5 to 20
\put {$\beta>0$} at 22 9
\setlinear
\setsolid
\plot 0 0 20 0 /
\plot 0 0 10 17.3 20 0 /
\setdashes
\plot 15 0 10 5.98 /
\plot 2.5 4.33 10 5.98 /
\plot 12.5 13 10 5.98 /
\setdots
\plot 10 5.98 6.14 10.60 /
\plot 10 5.98 15.80 7.26 /
\plot 10 5.98 7.88 0 /
\setsolid
\multiput {*} at 11.4 9.8 12.4 3 5.9 4.9 /
\multiput {$\bullet$} at 15 0 2.5 4.33 12.5 13 /
\put {$f_{12}^A$} at 14.0 3
\put {$f_{13}^A$} at 10 9.8
\put {$f_{23}^A$} at 5.9 6.5
\put {$R_{12}^A$} at 15 -1
\put {$R_{13}^A$} at 14.5 13
\put {$R_{23}^A$} at 1 4.33
\put {$Q^A_{12}$} at 5 10.60
\put {$Q^A_{23}$} at 16.80 7.26
\put {$Q^A_{13}$} at 7.88 -1
\put {$P^A_1$} at -1 0
\put {$P^A_3$} at 21 0
\put {$P^A_2$} at 10 18.3
\put {$2$} at  3 2
\put {$3$} at 10 14
\put {$1$} at 17 2
\arrow <5pt>  [0.2,0.4] from 10 5.5 to 12.0 3
\arrow <5pt>  [0.2,0.4] from 11.5 3 to 10.0 5.0
\setcoordinatesystem units <\dimen0,\dimen0> point at -28 0
\setplotarea x from 0 to 30, y from -5 to 20
\setsolid
\plot 0 0 20 0 /
\plot 0 0 10 17.3 20 0 /
\setdashes
\plot 15 0 10 5.98 /
\plot 2.5 4.33 10 5.98 /
\plot 12.5 13 10 5.98 /
\setdots
\plot 10 5.98 6.14 10.60 /
\plot 10 5.98 15.80 7.26 /
\plot 10 5.98 7.88 0 /
\setsolid
\multiput {*} at 11.4 9.8 12.4 3 5.9 4.9 /
\multiput {$\bullet$} at 15 0 2.5 4.33 12.5 13 /
\put {$f_{13}^B$} at 14.0 3
\put {$f_{23}^B$} at 10 9.8
\put {$f_{12}^B$} at 5.9 6.5
\put {$R_{13}^B$} at 15 -1
\put {$R_{23}^B$} at 14.5 13
\put {$R_{12}^B$} at 1 4.33
\put {$Q^B_{13}$} at 5 10.60
\put {$Q^B_{12}$} at 16.80 7.26
\put {$Q^B_{23}$} at 7.88 -1
\put {$P^B_1$} at -1 0
\put {$P^B_3$} at 21 0
\put {$P^B_2$} at 10 18.3
\put {$1$} at  3 2
\put {$2$} at 10 14
\put {$3$} at 17 2
\arrow <5pt>  [0.2,0.4] from 12.0 3 to  10 5.5
\arrow <5pt>  [0.2,0.4] from 10 5.5 to  5.9 4.4
\endpicture
\caption{\label{fig:ssb-bother}
{\small The simplices $\Sigma_A$ and $\Sigma_B$.  Various points on the 
sets $Z^A$ and $Z^B$ (and their extensions) are drawn.
$R^A_{23}=\frac{1}{2+\beta} (1+\beta,1,0),$
$Q^A_{12}=\frac{1}{1+2\beta} (\beta,1+\beta,0)$
(and similarly formulas for the other points
by cyclically permuting coordinates).
Similarly
$R^B_{12}=\frac{1}{2-\beta}(1,1-\beta,0)$
and
$Q^B_{13}=\frac{1}{1+\beta}(\beta,1,0).$
The points $f^A_{ij}$, $f^B_{kl}$ give the maximum extent
of the periodic orbit on $J$ ($\beta>\sigma$ only).
Take $(p^A,p^B)\in (Z^B_{1,2}\times Z^A_{3,1})$
as starting points.
Then $p^A(t)$
moves along $Z^B_{1,2}$ towards
$R^A_{12}$
and $p^B(t)$ moves
towards the midpoint of $\Sigma_B$ along $Z^A_{1,3}$.
Once $p^B(t)$ hits the midpoint of $\Sigma_B$ (provided $p^A(t)\ne
(1/3,1/3,1/3)$) then the flow extends continuously
and $(p^A(t),p^B(t))$ enters
$Z^B_{1,2}\times Z^A_{1,2}$.
Along this set $p^A(t)$ then moves back along $Z^B_{1,2}$
towards the midpoint in $\Sigma_A$ while $p^B(t)$ moves
towards $R^B_{12}$.
Next it enters the third part of $J$, and in this way
the orbit spirals (cyclically) around the different legs of $J$.
}}
\end{figure}

We shall prove here that trajectories on
$J$ visit the six pieces cyclically (in the order in which they appear 
in the definition in equation (\ref{eqn:defJ})), 
and the details of the behaviour 
are given in the following proposition (see Figure \ref{fig:ssb-bother}).

\begin{theo}[A third periodic which is created by a Hopf-like bifurcation]
\label{prop:onJ}
The system has a Hopf-like bifurcation when $\beta=\sigma$
on the topological two-dimensional manifold $J$.
Trajectories on $J$ spiral towards $(E^A, E^B)$
when $\beta\in (0,\sigma]$ while for
$\beta\in (\sigma,1)$ all trajectories in $J$
(except the one stationary at $(E^A, E^B)$)
tend to a periodic orbit $\Gamma$ in $J$.

We denote the endpoints of $\Gamma$ along the
pieces of $J$ by $f^A_{ij}$ and $f^B_{kl}$, as in Figure \ref{fig:ssb-bother}.
Then for $\beta\in (\sigma,1)$ one has the following formulas
for the `diameter' of the periodic orbit $\Gamma$:
$$\frac{|E^B_{kl}-f^B_{kl}|}{|Q^B_{kl}-f^B_{kl}|}=
\frac{|f^A_{ij}-m^A_{ij}|}{|R^A_{ij}-E^A_{ij}|}=
\frac{\beta^2+\beta-1}{2\beta+1}$$
and
$$\frac{|f^B_{kl}-E^B_{kl}|}{|R^B_{kl}-E^B_{kl}|}=
\frac{|E^A_{ij}-f^A_{ij}|}{|Q^A_{ij}-f^A_{ij}|}=
\frac{\beta^2+\beta-1}{(1+\beta)^2}.$$
Note that these expressions are both increasing
in $\beta\in (\sigma,1)$ and are zero when 
$\beta=\sigma$ (but unlike the usual Hopf-bifurcation
the diameter increases roughly linearly  when $\beta-\sigma>0$ 
is small).  They achieve their maximum $1/3$, $1/4$ both at $\beta=1$.
\end{theo}

The proof of this Theorem is a tedious calculation and
can be found in the appendix.
Note that it is also argued there that trajectories spiralling into 
$(E^A,E^B)$ on $J$ for $\beta<\sigma$, or
out of it for $\beta>\sigma$, do so in finite time. There is genuine
non-uniqueness of the flow at $(E^A, E^B)$ reflecting the fact
that it is the only point (for $\beta>0$) at which we cannot
continuously extend the flow.

\subsection{Dithering and chaotic behaviour near $J$}
\label{subsec:chaos}
It is not true that this orbit is attracting
for orbits that start off $J$. Nonetheless, orbits of the system
for $\beta > \sigma$ can spend long times near to this orbit $J$. 
The reason for this is that orbits can spiral along cones whose
axes consist of the 6 legs of the periodic orbit on $J$
in the way described in Figure~\ref{fig:quad}.
The next two figure below shows 
a portion of a numerically computed trajectory for $\beta=0.8$.
The orbit shown, starts near the centre, and during the time
interval considered comes close to the periodic orbit $J$,
while jittering a large number of times near this orbit.
Notice that the short time intervals and the large number of
time the players switch strategies.

\includegraphics[width=4in]{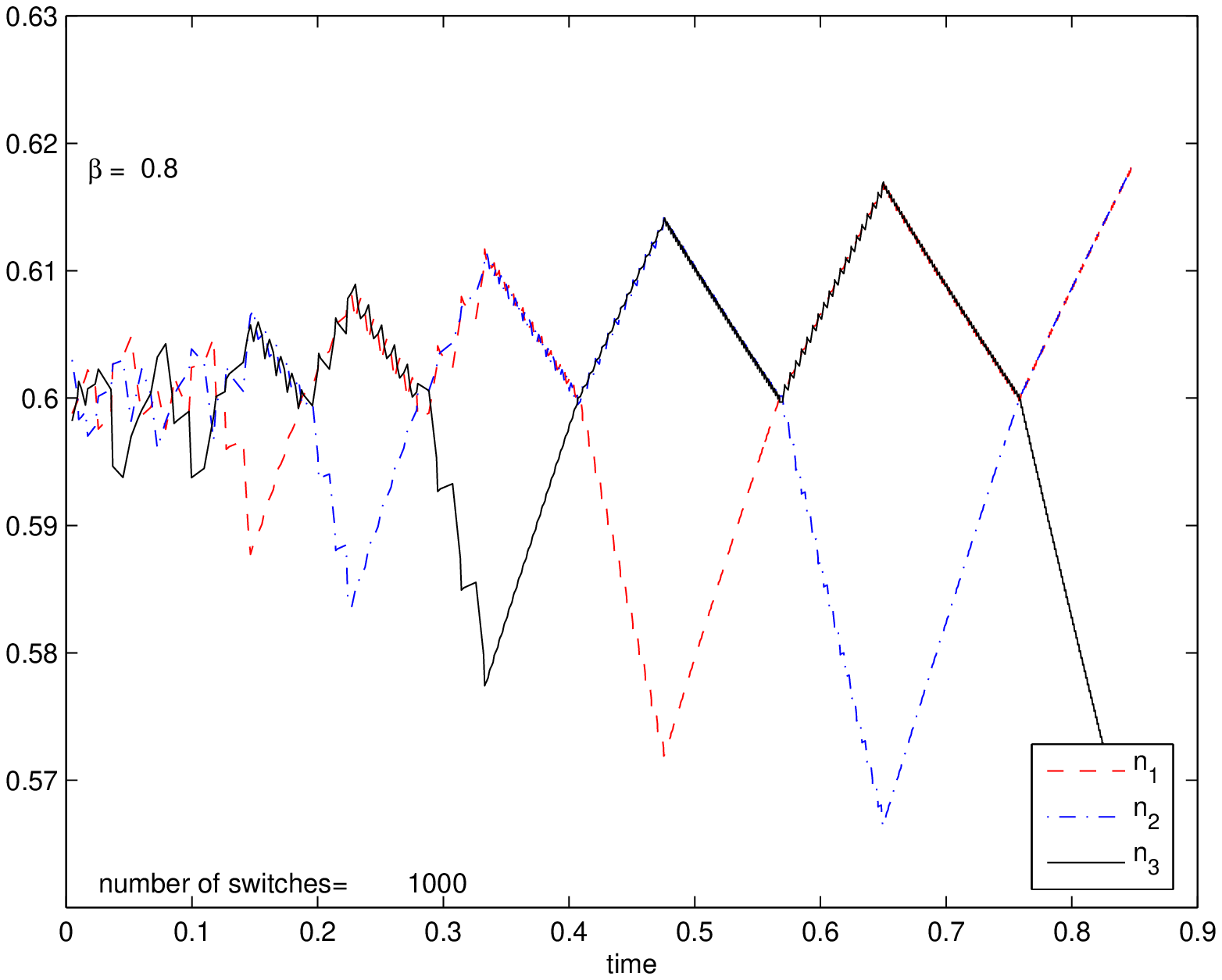}

\includegraphics[width=4in]{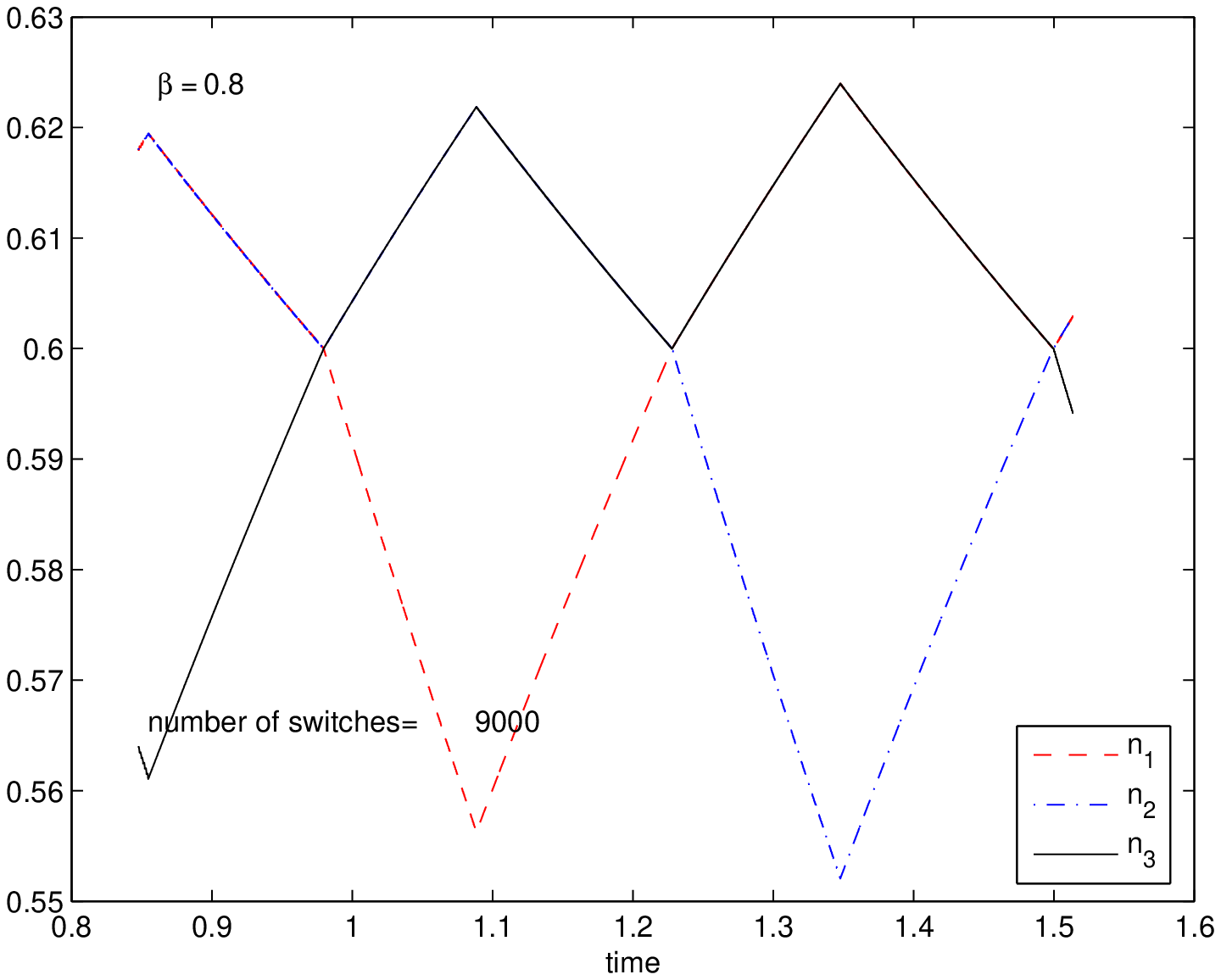}

The same orbit shown for a later time interval,
shows that orbits eventually become less periodic. 

\includegraphics[width=4in]{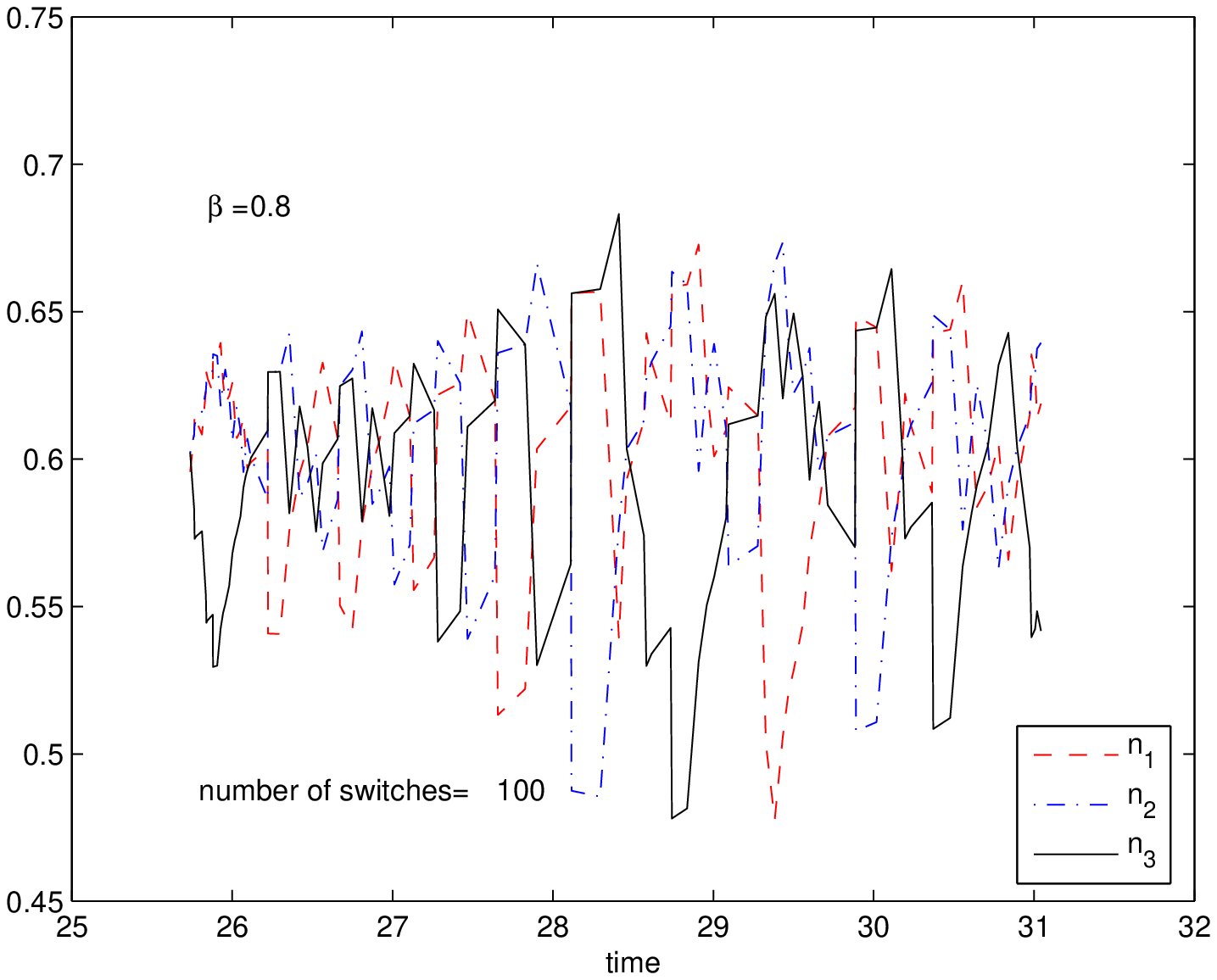}

In fact, as will be made precise in \citeasnoun{CS2006},
one has chaotic motion for $\beta\in [\sigma,\tau]$.
The figures below show orbits for $\beta=0.9$, both showing
erratic time-series and the attracting set projected
on $\Sigma^A$ and $\Sigma^B$.
%

\includegraphics[width=4in]{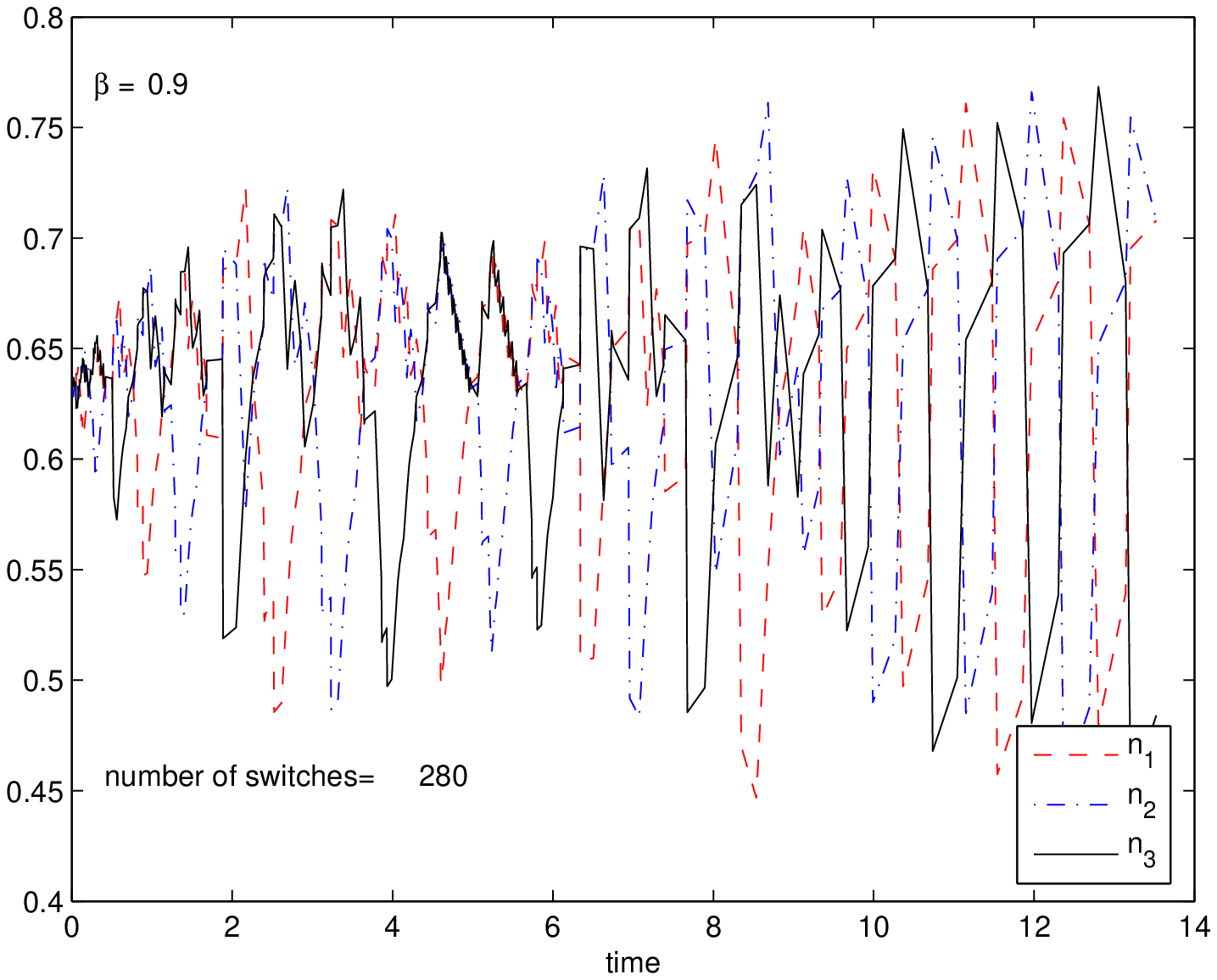}

\includegraphics[width=5in,height=4in]{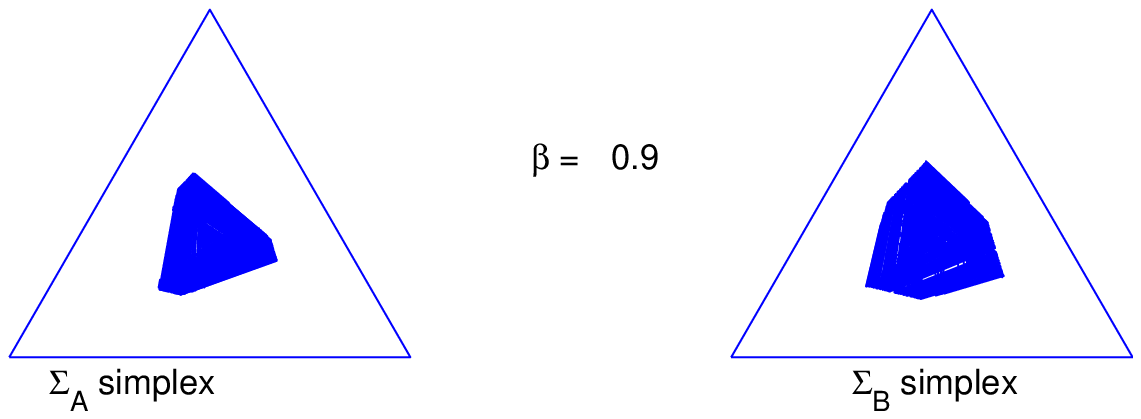}

\vskip -6cm

\section{Conclusion}
\label{s:concl}

For $\beta\in (-1,0]$ players always asymptotically 
become periodic.  When $\beta\in (0,\sigma)$
the Shapley orbit is still attracting but not globally attracting: 
there are orbits as described in 
Theorem~\ref{prop:onJ} which tend to $E$. In a sequel to this paper we shall show that there are infinitely many such cone-like sets
consisting of points which tend to $E$.
(So the stable manifold of $E$ is extremely complicated.)

The numerical results showed that for $\beta\in (\sigma,\tau)$ 
players act erratically and   chaotically.
In the sequel to this paper we shall give rigorous results
which explain this behaviour. The mathematical 
methodology in this sequel paper will be more geometric
(and will be based on first projecting the flow on 
a three-sphere and then by a study of the first return dynamics
to a global two-dimensional section consisting of a disc
bounded by the periodic orbit $\Gamma$).

Finally, when $\beta\in (\tau,1)$ orbits tend to a so-called
anti-Shapley orbit.

\bigskip

Of course most of the analysis we described also could be
applied to other families (with less symmetry). We will go into
this in the sequel to this paper.
Finally,  simulations and preliminary theoretical considerations indicate
that for $n \times n$ games with $n>3$ there are other routes
to chaotic behaviour than those studied here. 


\medskip

\noindent
{\bf Acknowledgements.} We are grateful to the anonymous
referee and an associate editor for helpful comments and suggestions.


\bibliography{games.bib}

\section{Appendix}
\label{appendix}

In this appendix we collect together a number of explicit calculations 
that are used in the main text. 

It is convenient to work in the space $(v^A, v^B)$ where $v^A = A p^B$ and 
$v^B = p^AB$ (so $v^A$ is a column vector, and $v^B$ is a row
vector). From the description of the dynamics in the main text,
when $v^A_i$ is the largest coordinate of $v^A$, 
player $A$ adjusts $p^A$ towards $P^A_i$, which implies $v^B$ moves in a 
straight line 
towards $B_i$, the $i^{\rm th}$ row of $B$. Similarly,
if $v^B_j$ is the largest component of $v^B$, $v^A$ moves towards $A_j$, the
$j^{\rm th}$ column of $A$. In both cases, our convention
is that the motion occurs at a speed 
such that from a given initial condition, and assuming that
$i$ and $j$ were to remain the maximal coordinates of $v^A$ and $v^B$, the
vectors would reach $B_i$ and $A_j$ in one unit of time. 
So, for $0\le s \le 1$,
provided the maximal coordinate of $v^A$ or $v^B$ does not change, we can
write
\begin{equation}
\begin{array}{rc}
 v^A(t+s) &= v^A(t) + s (A_j - v^A(t)) \\
 v^B(t+s) &= v^B(t) + s (B_i - v^B(t))
\end{array} \,.
\label{eqn:app1}
\end{equation}
Each time that
the maximal coordinate of either $v^A$ or $v^B$ changes, we reset $s$ in
(\ref{eqn:app1}) to
zero.  In this appendix we calculate explicit orbits with particular
properties, using these equations. We illustrate our 
calculations with figures that show the development of the components 
of $v^A$ and $v^B$ with time. 
Let $\sigma=(-1+\sqrt{5})/2$, the golden mean.

\subsection{Existence of the clockwise Shapley orbit in $\beta \in (-1, \sigma)$}
\label{ss:appclock}

\begin{figure}[htp]
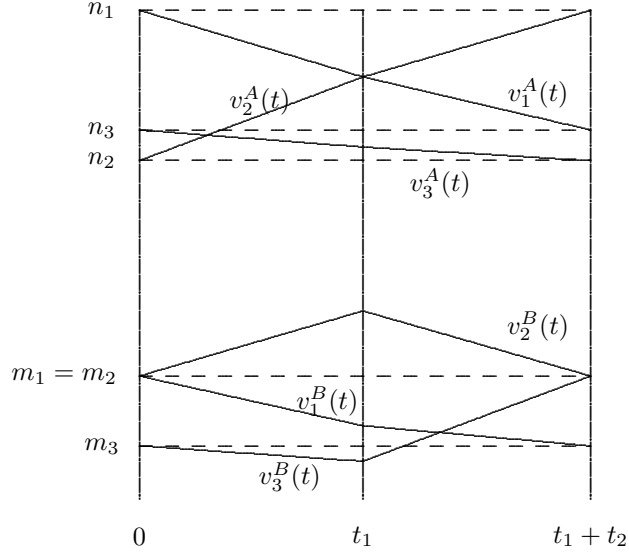
 \hfil
\beginpicture
\dimen0=0.1cm
\setcoordinatesystem units <\dimen0,\dimen0>
\setplotarea x from 0 to 60, y from 0 to 70
\setlinear

\plot 0 5 0 70 /
\plot 29.72 5 29.72 70 /
\plot 60 5 60 70 /
\put {$0$} at 0 0 \put{$t_1$} at 29.72 0  \put {$t_1+t_2$} at 60 0


\plot 0 21.338 29.72 14.72 60 12.044 / 
                  \put{$m_1=m_2$} at -10 21.338 \put {$v^B_1(t)$} at  25 18
\plot 0 21.338 29.72 30 60 21.338 / 
                  \put {$v^B_2(t)$} at 53 28
\plot 0 12.044 29.72 10 60 21.338 / 
                  \put {$v^B_3(t)$} at 20 8
                  \put {$m_3$} at -5 12.044


\plot 0 70 29.72 61.130  60 54.07     /  
          \put {$n_1$} at -5 70 \put {$v^A_1(t)$} at 53 59
\plot 0  50 29.72 61.130 60  70 /  
          \put {$n_2$} at -5 50 \put {$v^A_2(t)$} at 16 58  
\plot 0  54.07 29.72 51.801 60  50    /  
          \put {$n_3$} at -5 54.07 \put {$v^A_3(t)$} at 40 47

\setdashes
\plot 0 50 60 50 / \plot 0 54.07 60 54.07 / \plot 0 70 60 70 /
\plot 0 21.338 60 21.338 / \plot 0 12.044 60 12.044 /

\endpicture
\caption{\label{fig:nmclock}
{\small Components of $v^A$ and $v^B$ plotted against time for one 
third of the symmetric Shapley (clockwise) periodic orbit 
which exists for
$\beta\in [0, \sigma)$. 
Since $v^B_2$ is maximal throughout, $v^A(t)$
moves towards $A_2 = (0,1,\beta)^{T}$. 
In the interval $[0,t_1]$
coordinate $v^A_1$ is maximal and $v^B(t)$ moves towards $B_1=(-\beta,1,0)$; 
in the interval $[t_1, t_1+t_2]$, $v^A_2$ is maximal and $v^B(t)$ 
moves towards $B_2= (0,-\beta,1)$.}}
\end{figure}

Let us compute the clockwise Shapley periodic orbit
which, as we will see, exists for $\beta\in (-1,\sigma)$.
The periodic orbit we are looking for is rotation symmetric, 
so we will only compute the part of the orbit 
where player B prefers strategy 2 (while player A
prefers first strategy 1 and then strategy 2); see Figure 
\ref{fig:nmclock}.

Let
$n=(n_1,n_2,n_3)^{T}$ and $m=(m_1,m_2,m_3)$ be the initial values
of $v^A$ and $v^B$ respectively. 
Then we require, as shown in the figure, that 
\begin{equation}
n_1>\max(n_2,n_3) \quad \mbox{ and } \quad m_1=m_2>m_3\,,
\label{eqn:app13}
\end{equation}
so that initially $v^A$ heads 
towards $A_2= (0,1,\beta)^{T}$
and $v^B$ towards $(-\beta,1,0)=B_1$.
The first time $t_1$ that the orbit meets $Z$ again
is when
the first and second coordinates
of $v^A(t)=(1-t)n^{T}+t(0, 1, \beta)^{T}$ become equal which occurs
when $t=t_1>0$,
$$t_1=\frac{n_1-n_2}{n_1-n_2+1}.$$
(Of course, it is possible that we first get that $n_1=n_3$, 
but this does not occur because $n_1=n_3$ happens
when $t'=\frac{n_1-n_3}{n_1-n_3+\beta}$ which is larger than
$t_1$ as $\beta<1$.)
The $v^A$ and $v^B$ values at $t=t_1$ are
$$v^A(t_1)=\frac{1}{(n_1-n_2+1)}\left(\begin{array}{c}
n_1\\ n_2+(n_1-n_2)\\ n_3+(n_1-n_2)\beta \end{array}\right) $$ 
and
$$v^B(t_1)=\frac{1}{(n_1-n_2+1)}\left(m_1-\beta(n_1-n_2),\, m_2+(n_1-n_2),
\, m_3\right)\,.$$
Next we have that $t_2>0$
is the smallest number so that the second and third coordinate
of $v^B(t_1+t_2)=(1-t_2)v^B(t_1)+t_2(0,-\beta,1)$ are equal. This gives 
$$t_2=\frac{m_2(t_1)-m_3(t_1)}{m_2(t_1)-m_3(t_1)+1+\beta}=
\frac{(m_2-m_3)+(n_1-n_2)}
{(m_2-m_3)+(n_1-n_2)+(1+\beta)(1+n_1-n_2)}$$
which implies
\begin{equation}
v^A(t_1+t_2)=\frac{1}{X}  \left( \begin{array}{c}
n_1(1+\beta) \\ n_1(1+\beta)+\delta \\
((n_3+(n_1-n_2)\beta)(1+\beta)+\beta\delta) \end{array}\right)\,,
\label{eqn:app14}
\end{equation}
and
\begin{equation}
v^B(t_1+t_2)=\frac{1+\beta}{X}\left((m_1-\beta(n_1-n_2)),\, 
(m_2 + (n_1-n_2))-\beta\delta,\,m_3+\delta/(1+\beta) \right)
\label{eqn:app15}
\end{equation}
where $\delta=(m_2-m_3)+(n_1-n_2)$
and $X=\delta+(1+\beta)(1+n_1-n_2)$. 
For a symmetric periodic orbit we have to satisfy 
\begin{equation}
v^A(t_1+t_2)=
\left(\begin{array}{c} n_3\\ n_1\\ n_2\end{array} \right) 
\quad \mbox{ and } \quad 
v^B(t_1+t_2) = 
\left(\begin{array}{c} m_3\\ m_1\\ m_1\end{array} \right)\,, 
\label{eqn:app16}
\end{equation}
and the strict inequalities (\ref{eqn:app13}).
Solving equations (\ref{eqn:app14})--(\ref{eqn:app16}) by hand
(or by a MAPLE worksheet \footnote{which can be downloaded from
http://www.maths.warwick.ac.uk/$\sim$strien/Publications} gives
\begin{equation}n_1 =\nu\label{eq:solclock1}\end{equation}
\begin{equation}n_2 =-\frac{
  - 2 \beta^2 - \beta^3 + 2\,\nu \beta^3 - 3\,\nu\, \beta - 2\,\nu 
+ 3\,\nu^2 + 3\,\nu^{2}\,\beta^2 + 3\,\nu^2\, \beta}{1+\beta+\beta^2},
\label{eq:solclock2}\end{equation}
\begin{equation}n_3 =\frac{2\, \beta + 1 + 2\,\nu\, \beta^3 - \nu\,\beta^2 
- 4\,\nu\,\beta - 3\,\nu + 3\,\nu^2 + 3\,\nu^2\,\beta^2 + 3\,\nu^2\,\beta}
{1+\beta+\beta^2},\label{eq:solclock3}
\end{equation}
\begin{equation}m_2 =m_1=1-\nu-\nu\beta \label{eq:solclock4}
\end{equation}
\begin{equation}m_3 =-\beta+2\nu\beta-1+2\nu,\label{eq:solclock5}
\end{equation}
where $\nu$ satisfies $f_\beta(\nu)=0$ with
\begin{equation}
f_\beta(Z)=\left(3\beta^2 + 3\beta + 3\right)  Z^3
 +  \left( 2\beta^3 - 2\beta^2 - 5\beta - 4 \right) Z^2
+ \left( -\beta^3+4\beta + 3\right) Z - 1-\beta\,.
\label{eqn:clockpoly}
\end{equation}
Note that for $\beta\in (-1,1)$, $f_\beta(0)=-1-\beta<0$ and
$f_\beta(1)=\beta^3+\beta^2+\beta+1>0$. 
For $\beta\in (-1,1)$, the function $t\mapsto f_\beta(t)$
is monotone because its derivative
\begin{equation}
f_\beta'(Z)=3\left(3\beta^2 + 3\beta + 3\right)  Z^2
 +  2\left( 2\beta^3 - 2\beta^2 - 5\beta - 4 \right) Z
+ \left( -\beta^3+4\beta + 3\right) 
\end{equation}
is a strictly postive quadratic function;
this is because  its discriminant is
equal to $16\beta^6+4\beta^5-28\beta^4-92\beta^3-88\beta^2-92\beta-44$
and so is negative when $\beta\in (-1,1)$.
Hence, for $\beta\in (-1,1)$ the
cubic (\ref{eqn:clockpoly})  has a unique real root $\nu$, and this root
is in $(0,1)$. Moreover, this solution $\nu$ depends
continuously on $\beta$. When
$\beta=0$, this solution satisfies the inequalities that guarantee the
behaviour is as shown in the figure: indeed when $\beta=0$ then
we obtain 
$n_1=0.594..$, $n_2=0.129..$, $n_3=0.277..$, $m_1=m_2=0.405..$,
$m_3=0.188..$. 

Let us show that equalities can only hold in
(\ref{eqn:app13})  when $\beta\notin (-1,1)$ or $\beta=\sigma$.
If $m_1=m_3$ then the expressions for
$m_1$ and $m_3$ give
$\nu_\beta=(2+\beta)/(3+3\beta)$ and then
$$f_\beta(\nu_\beta)=-\frac{(\beta^2+\beta-1)(1+\beta+\beta^2)^2}{9(1+\beta)^3}=0;$$
this implies that either $\beta\notin (-1,1)$ or $\beta$ is equal to the golden mean $\sigma$.
If $n_1=n_3$ then because of the first equation
in (\ref{eqn:app16}) the first two components $v^A(t_1+t_2)$ in
(\ref{eqn:app14}) are equal, and so $\delta=0$; using these two equation again we get 
$n_1-n_2=0$ and with $\delta=0$ this implies $m_1=m_3$
and therefore as we saw $\beta\notin (-1,1)$ or $\beta=\sigma$. 
If $n_1=n_2$ then 
$X=\delta+1+\beta$ and
$$
v^A(t_1+t_2)=\frac{1}{X}  \left( \begin{array}{c}
n_1(1+\beta) \\ n_1(1+\beta)+\delta \\
n_3(1+\beta)+\beta\delta\end{array}\right) = \left( \begin{array}{c}
n_3 \\ n_1 \\ n_1 \end{array} \right)
$$
which gives that $\delta n_3 + \delta=\beta \delta$
and so either $\delta=0$ or $n_3=\beta-1$ and so $\nu=n_1=n_2=1$.
The latter is impossible because $f_\beta(1)>0$ while the latter implies that $m_2=m_3$ ( we assumed that $n_1=n_2$);
 as we showed above $m_2=m_3$
implies $\beta\notin (-1,1)$ or $\beta=\sigma$.

It follows that there exists a symmetric clockwise periodic orbit
for all $\beta\in (-1,\sigma)$, depending continuously on $\beta$. 
As $\beta\uparrow \sigma$
this periodic orbit shrinks to the equilibrium solution
$n_1=n_2=n_3=(1+\beta)/3$ and $m_1=m_2=m_3=(1-\beta)/3$.
To see this, notice that
$f_\beta((1+\beta)/3)=0$ is
a fifth order equation which can be factorized as
$$\frac{1}{3}(1+\beta)(\beta^2+\beta-1)(\beta^2-\beta+1)=0$$
which is zero if $\beta=\sigma$ or $\beta\notin (-1,1)$.
If $\beta\in (\sigma,1)$ the inequalities in (\ref{eqn:app13}) no longer
hold. This can  be seen by considering the solutions for $\beta=1$ 
in (\ref{eq:solclock1}-\ref{eq:solclock5}). 

\subsection{Existence of the anticlockwise 
periodic orbit when $\beta\in (\sigma,1)$.}
\label{ss:appanti}

We first show that for 
$\sigma  < \beta \le 1$
there exists a unique symmetric anticlockwise periodic 
orbit as described in Theorem~\ref{exipe}in the main text. 

For such an orbit to exist, with the appropriate symmetry,
the vectors $v^A$ and  $v^B$ must behave as shown in
Figure \ref{fig:nmanti} below. 

\begin{figure}[htp]
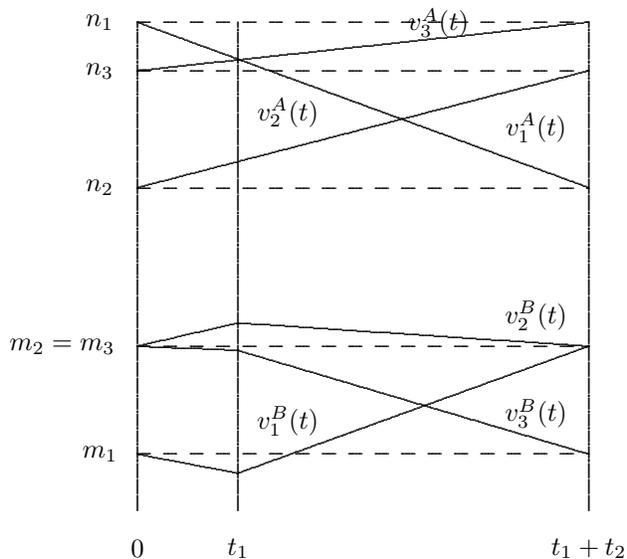
 \hfil
\beginpicture
\dimen0=0.1cm
\setcoordinatesystem units <\dimen0,\dimen0>
\setplotarea x from 0 to 60, y from 0 to 70
\setlinear

\plot 0 5 0 70 /
\plot 13.36 5 13.36 70 /
\plot 60 5 60 70 /
\put {$0$} at 0 0 \put{$t_1$} at 13.36 0  \put {$t_1+t_2$} at 60 0


\plot 0 12.542 13.36 10 60 26.906 / 
                  \put{$m_1$} at -5 12.542 \put {$v^B_1(t)$} at  20 17
\plot 0 26.906 13.36 30 60 26.906 / 
                  \put {$m_2=m_3$} at -10 26.906
                  \put {$v^B_2(t)$} at 53 31
\plot 0 26.906 13.36 26.34 60 12.542 / 
                  \put {$v^B_3(t)$} at 53 18


\plot 0 70       60 48     /  
          \put {$n_1$} at -5 70 \put {$v^A_1(t)$} at 53 56
\plot 0 48       60 63.554 /  
          \put {$n_2$} at -5 48 \put {$v^A_2(t)$} at 20 58  
\plot 0 63.554   60 70     /  
          \put {$n_3$} at -5 63.554 \put {$v^A_3(t)$} at 40 70

\setdashes
\plot 0 48 60 48 / \plot 0 63.554 60 63.554 / \plot 0 70 60 70 /
\plot 0 26.906 60 26.906 / \plot 0 12.542 60 12.542 /

\endpicture
\caption{\label{fig:nmanti}
{\small Components of $v^A$ and $v^B$ plotted against time for one 
third of the symmetric (anticlockwise) periodic orbit which exists for
$\beta\in (\sigma,1)$. Note that
the identity of the maximal coordinate of $v^A$ determines the motion of 
$v^B$, and vice versa.  Since $v^B_2$ is maximal throughout, $v^A(t)$
moves towards $A_2 = (0,1,\beta)^{T}$ in both time intervals. 
In the interval $[0,t_1]$
coordinate $v^A_1$ is maximal and $v^B(t)$ moves towards $B_1=(-\beta,1,0)$; 
in the interval $[t_1, t_1+t_2]$, $v^A_3$ is maximal and $v^B(t)$ 
moves towards $B_3= (1,0,-\beta)$. }}
\end{figure}

We compute only the part of the orbit where player B prefers strategy $2$,
starting where player B is
indifferent between strategies $2$ and $3$, and ending where player
$B$ is indifferent between strategies $1$ and $2$. During this interval,
player A prefers first strategy $1$ and then strategy $3$.  
We then use the
rotational symmetry to determine whether this corresponds to a
periodic orbit. The complete orbit consists of three such sections, 
with the coordinates permuted cyclically on each section.

As before, to avoid too many sub and superscripts, let
$n=(n_1,n_2,n_3)^{\rm T}$ be the initial value  $v^A(0)$,
and $m=(m_1,m_2,m_3)$  be the initial value  $v^B(0)$.
So the starting position given in the
figure satisfies 
\begin{equation}
n_1>n_2,\, n_3 \quad \mbox{ and } \quad m_2=m_3>m_1\,,
\label{eqn:app2}
\end{equation}
and the orbit in $v^A, v^B$ space initially
heads to respectively $(0,1,\beta)^{T}=A_2$
and $(-\beta,1,0)=B_1$ as we require.
The first time that the orbit meets $Z$ (the indifference
set) again in $t>0$ is when
the first and third coordinates of 
\begin{equation}
v^A(t)=(1-t)n+t\left(\begin{array}{c} 0\\ 1\\ \beta\end{array} \right) 
\label{eqn:app3}
\end{equation}
become equal and player A becomes indifferent between
strategies $1$ and $3$. Solving $v^A_1(t)=v^A_3(t)$ in
(\ref{eqn:app3}) for $t=t_1>0$,
we obtain
\begin{equation}
t_1=\frac{n_1-n_3}{n_1-n_3+\beta}.
\label{eqn:app4}
\end{equation}
Of course, depending on the values of the $n_i$ it is possible 
that we first get that $n_1(t)=n_2(t)$
for some $t\in (0,t_1)$.
Since we are not interested in this case, we insist that this does
not occur and so our initial 
condition satisfies
\begin{equation}
(n_1-n_2)\beta>(n_1-n_3). 
\label{eqn:app5}
\end{equation}
Note that the motion of $v^B$ ensures that
the largest coordinate of $v^B$ must be $v^B_2$ for the whole of the period
$(0,t_1)$, so there is no danger that player B becomes indifferent before
player A.

The $v^A$ and $v^B$ values at $t=t_1$ are, using (\ref{eqn:app1}),
\begin{equation}
v^A(t_1)=\frac{1}{(n_1-n_3+\beta)}\left(\begin{array}{c}
\beta n_1 \\ \beta n_2 + n_1-n_3 \\ \beta n_1\end{array} \right)
\label{eqn:app6}
\end{equation}
and
\begin{equation}
v^B(t_1)=\frac{1}{(n_1-n_3+\beta)}\left(\beta m_1-\beta(n_1-n_3),\ \ 
\beta m_2 + (n_1-n_3),\ \ \beta m_3\right).
\label{eqn:app7}
\end{equation}
After this moment, player $A$ changes direction
and so $v^B=p^AB$ heads for $(1,0,-\beta)$
and  $t_2>0$ becomes
the smallest number so that the first and second coordinates
of $v^B(t_1+t_2)=(1-t_2)v^B(t_1)+t_2(1,0,-\beta)$ become equal. 
This occurs when
\begin{equation}
t_2=\frac{v^B_2(t_1)-v^B_1(t_1)}{v^B_2(t_1)-v^B_1(t_1)+1}=
\frac{\beta(m_2-m_1)+(1+\beta)(n_1-n_3)}
{\beta(m_2-m_1)+(2+\beta)(n_1-n_3)+\beta}.
\label{eqn:app8}
\end{equation}
Hence, writing $$\delta=\beta(m_2-m_1)+(1+\beta)(n_1-n_3)$$
we have
$$t_2=\frac{\delta}{\delta+\beta+n_1-n_3} \,.$$
Using (\ref{eqn:app1}) again, this implies
\begin{equation}
v^A(t_1+t_2)=\frac{1}{\delta+\beta+n_1-n_3}\left(
\begin{array}{c}
n_1\beta \\ n_1-n_3+n_2\beta+\delta \\ n_1\beta+\beta\delta\end{array}
\right)\,,
\label{eqn:app9}
\end{equation}
and
\begin{equation}
v^B(t_1+t_2)=\frac{1}{\delta+\beta+n_1-n_3}\left((\beta m_2+n_1-n_3),\,
(\beta m_2 + n_1-n_3),\, \beta(m_2-\delta) \right)\,.
\label{eqn:app10}
\end{equation}
To get an anticlockwise symmetric periodic orbit we require that
\begin{equation}
v^A(t_1+t_2)=
\left(\begin{array}{c} n_2\\ n_3\\ n_1\end{array} \right) 
\quad \mbox{ and } \quad 
v^B(t_1+t_2) = 
\left(\begin{array}{c} m_2\\ m_2\\ m_1\end{array} \right)\,. 
\label{eqn:app11}
\end{equation}
Combining (\ref{eqn:app9}), (\ref{eqn:app10}) and (\ref{eqn:app11}),
we need to solve the following equations:
\begin{equation}
\begin{array}{rc}
n_1\cdot \beta &= n_2\cdot (\delta+n_1-n_3+\beta)\\
n_1-n_3+n_2\cdot \beta+\delta&=n_3\cdot (\delta+n_1-n_3+\beta)\\
n_1\cdot \beta+\beta\cdot \delta&=n_1\cdot (\delta+n_1-n_3+\beta)\\
m_2\cdot \beta+n_1-n_3&=m_2\cdot (\delta+n_1-n_3+\beta)\\
m_2\cdot \beta-\beta\cdot \delta&=m_1\cdot (\delta+n_1-n_3+\beta)
\end{array}
\label{eqn:anti}
\end{equation}
and, collecting together (\ref{eqn:app2}) and (\ref{eqn:app5}), 
we require solutions to (\ref{eqn:anti}) that satisfy
the strict inequalities 
\begin{equation}
m_2=m_3>m_1, \qquad n_1>\max(n_2,n_3)
\qquad {\rm and}\qquad (n_1-n_2)\beta>n_1-n_3\,.
\label{eqn:app12}
\end{equation}

It is straightforward to check using the computer algebra package 
MAPLE,\footnote{which can be downloaded from
http://www.maths.warwick.ac.uk/$\sim$strien/Publications} or
more laboriously by hand, 
\footnote{
If manipulating these equations by hand, it is useful to observe 
that the solution must satisfy
$$n_1+n_2+n_3 = 1 + \beta \qquad \mbox{ and } \qquad m_1+m_2+m_3=1-\beta$$
which come from the definitions of $v^A$ and $v^B$ and the fact 
that $p^A$ and $p^B$ are probability vectors.}
that solutions to the equalities 
(\ref{eqn:anti}) satisfy:
$$m_2=m_3=\mu,$$
$$m_1= -  \beta + 1 - 2\,\mu,$$
$$n_1=\beta-\mu \beta$$
$$n_2=
\frac{2\,\beta^2 + 1 - 3\,\mu\, \beta^3 - 5\,\mu \beta^2 - 2\,\mu\,\beta -
2\,\mu + 3\,\mu^2\, \beta^2 + 3\,\mu^2\, \beta + 3\,\mu^2\,\beta^{3} }{
\beta + 1 +  \beta^2}.$$
$$n_3=- \frac{\beta^2 -  \beta - 4\,\mu\, \beta^3
  - 6\,\mu\, \beta^{2} - 3\,\mu\,\beta - 2\,\mu + 3\,\mu^2\,\beta^2 
+ 3\,\mu^2\,\beta + 3\,\mu^2\,\beta^3 }{\beta+ 1 + \beta^{2}} ,$$
where $\mu$ satisfies $f_\beta(\mu)=0$ and where
$$ f_\beta(Z)=\left(3\beta^2 + 3\,\beta + 3\,\beta^3\right) \, Z^3
 +  \left( - 5\,\beta^3 - 7\,\beta^2 - 4\,\beta - 2\, \right) \,Z^2+$$
\begin{equation}
+ \left( 1 + \beta + 5\beta^2 + 2\beta^{3}\right) Z - { \beta}^{2}\,.
\label{eqn:antipoly}
\end{equation}
Note that for $\beta\in (0,1)$, $f_\beta(0)=-\beta^2<0$ and $f_\beta(1)=-2<0$
while $f_\beta(1/3)=(2/9)\beta^3+1/9>0$. This implies that the
cubic polynomial (\ref{eqn:antipoly}) always has three roots
$\mu$,  two in $(0,1)$ and one in $(1,\infty)$. When
$\beta=1$, the only of these three roots which leads to a 
solution satisfying the inequalities
(\ref{eqn:app12}) has $\mu=0.155..$, so
$n_1=0.844..$, $n_2=0.449..$, $n_3=0.706..$, $m_1=-0.311..$, and
$m_2=m_3= 0.155..$. This solution (which depends on $\beta$) continues 
to exist as $\beta$ decreases, and gives a periodic orbit of the
desired type until one of the strict inequalities (\ref{eqn:app12})
fails to hold. 
Let us check what happens if one of these inequalities fails, still assuming $\beta\in (0,1)$.
If $m_1=m_2$ then by the expressions above, $\mu=(1-\beta)/3$.
If $n_1=n_2$ then by the first and third 
equations in (\ref{eqn:anti}), $\delta=0$ and either $n_1=n_2=n_3$ or 
$n_1=0$ and so $\mu=1$ (which is impossible because $1$ is not a root of $f_\beta$).
If $n_1=n_3$ then the first three equations in (\ref{eqn:anti}) imply 
$n_1=\beta$ or $\delta=0$;
if $n_1=\beta$ then using the formula for $n_1$ we get
$\beta=\beta-\mu\beta$ and so either $\beta=0$ or $\mu=0$ (both of which 
are impossible) whereas  if $\delta=0$ (and $n_1=n_3$) 
then $m_1=m_2$ and so, from the previous argument, $\mu=(1-\beta)/3$.
If $(n_1-n_3)\beta=n_1-n_2$ then manipulating the first three equations
in (\ref{eqn:anti}) we get $1/(1-\beta)=\beta$ (which is impossible
since $\beta\in (0,1)$) or  $\delta=0$ which in turn implies $n_1=n_2=n_3$ 
or $n_1-n_2+1=0$ which (after some manipulations) gives $n_1=0$ and so
$\mu=0$ which again is impossible.
In other words, the only way the inequalities can break down is
when $\mu$ becomes equal to $\mu=(1-\beta)/3$.
But if we substitute $Z=(1-\beta)/3$ into (\ref{eqn:antipoly}) we obtain
a sixth order equation which can be factorized as
$$\frac{-1}{9}(\beta_c^2+\beta_c-1)(\beta_c+1+\beta_c^2)^2 = 0\,.$$
Hence the critical value $\beta_c$ of $\beta$ is equal to the golden mean $\sigma$ 
(since the other solutions of
the sixth order polynomial in $\beta$ are negative or complex). It is
straightforward to check that for all $\beta \in (0,1)$ the 
solution generated by the other root of (\ref{eqn:antipoly}) fails to satisfy
the strict inequalities (\ref{eqn:app12}). 
Thus we have shown that there is unique
symmetric anticlockwise periodic orbit as described in the proposition
for $\beta \in (\sigma, 1]$ and none
otherwise. Furthermore, the orbit shrinks down to the interior
equilibrium as $\beta \downarrow \sigma$.

\subsection{Stability of the anticlockwise and clockwise periodic orbits}
\label{ss:appstab}

As stated in the main text:

\begin{prop}
Shapley's periodic orbit is attracting for all
$\beta\in [0,\sigma)$. The anticlockwise orbit which exists
for $\beta\in (\sigma,1)$ is attracting when
$\beta\in (\tau,1)$ and of saddle-type when
$\beta\in (\sigma,\tau)$. Here $\tau$ is a root of a specific expression whose value is a little larger than $0.915$.
\end{prop}
\begin{proof}
Let us just show  the computation for the anticlockwise periodic orbit, because
in that case an interesting bifurcation occurs.
Since $\sum n_i=1+\beta$ and $\sum m_i=1-\beta$,
we shall perturb the initial conditions
considered in Section \ref{ss:appanti} above to
$n_1'=n_1+\epsilon_1$, $n_2'=n_2+\epsilon_2$, $n_3'=n_3-\epsilon_1-\epsilon_2$
while
$m_1'=m_1-2\epsilon_3$, $m_2'=m_3'=m_2+\epsilon_3$.
Repeating the calculation of Section \ref{ss:appanti}
with these perturbed values we find that 
after slightly altered times
$t_1'$ and $t_2'$ we arrive at a final position $v^A(t_1'+t_2')$ 
with coordinates
$v^A_1(t_1+t_2)+\gamma_1$,
$v^A_2(t_1+t_2)+\gamma_2$, $v^A_3(t_1+t_2) -\gamma_1-\gamma_2$
and a final position for $v^B(t_1'+t_2')$ with coordinates
$v^B_1(t_1+t_2)+\gamma_3$, $v^B_1(t_1+t_2)+\gamma_3$, 
$v^B_3(t_1+t_2)-2\gamma_3$
where
$$\left(
\begin{array}{c}
\gamma_1 \\ \gamma_2 \\ \gamma_3 \end{array}\right)= 
\frac{n_2}{n_1\beta}
\left( \begin{array}{ccc}
\beta(3+2\beta)-n_1 2(2+\beta) & \beta(1-n_1+\beta)-2n_1
& 3\beta^2-3n_1\beta \\
\quad \beta-2n_2 (2+\beta) \quad & -n_2(2+\beta) & -3n_2 \beta \\
2-\frac{2(\beta-n_1)(2+\beta)}{\beta} & 
\quad 1-\frac{(\beta-n_1)(2+\beta)}{\beta} \quad
& 3 n_1 - 2\beta
\end{array}\right)
\left(\begin{array}{c}
\epsilon_1 \\ \epsilon_2 \\ \epsilon_3 \end{array}\right)
$$
plus higher order terms. 
\footnote{This matrix can be deduced from equations (\ref{eqn:anti}),
where it is convenient to observe that we can use
the first of those equations to write
$(\delta+n_1-n_3+\beta) = \beta n_1 / n_2$.} 
One can compute the eigenvalues of this matrix
using MAPLE; one
of of them is $n_2/n_1$,
while the others are given as rather complicated (but exact) formulae
in $n_1$, $n_2$ and $\beta$: 
\begin{equation}
\frac
{n_2
\left(-2n_1\beta-n_1-2n_2-n_2\beta+2\beta^2+\pm \sqrt{\Delta}\right)}
{2n_1\beta},
\label{eq:eig23}
\end{equation}
where
\begin{equation*}
\begin{array}{rl} \Delta &:=10\,n_1\,n_2\,\beta+4\,n_1\,n_2-4\,n_2\,{\beta}^{2}-4\,n_2\,{\beta}^{3}-8\,n_1\,{\beta}^{3}+4\,{\beta}^{4}+4\,{n_2}^{2}+4\,{n_2}^{2}\beta+\\
&\quad 4\,{n_1}^{2}{\beta}^{2}+4\,{n_1}^{2}\beta+{n_2}^{2}{\beta}^{2}
\mbox{}+4\,n_2\,{\beta}^{2}n_1-8\,n_1\,{\beta}^{2}+{n_1}^{2}-4\,n_2\,\beta-8\,n_1\,\beta+4\,{\beta}^{2}+4\,{\beta}^{3}
\end{array}
\end{equation*}
Using the explicit formulae for
$n_i$ computed in Section \ref{ss:appanti}, and by inserting these 
in the expressions for the eigenvalues,
we find that the eigenvalues are expressions in $\beta$. 
When $\beta=1$, these are 
$0.532..$, $-0.815..$, and $-0.184..$ (where the first one
is $n_2/n_1$). One can check  that as $\beta$ varies
the eigenvalues (\ref{eq:eig23}) always remain negative.
By numerically evaluating the expressions for these (which
can be done as accurately as required) we find that
one of them crosses $-1$ when $\beta$ is equal to $\tau$
where $\tau$ is roughly $0.915$
while the other remains, for all $\beta\in (\sigma,1)$,
in $(-1,0)$.
Consideration of the sizes and signs of the
eigenvalues described shows
that for $\beta>\tau$ the periodic orbit is attracting,
and for $\beta<\tau$ it becomes of saddle-type.

Note that we do {\bf not} have a generic periodic doubling
bifurcation 
\footnote{as would be expected in a normal differential
equation}
when one of the eigenvalues becomes equal to $-1$.
Indeed, at this parameter the linear part of
the above transformation has an eigenvector $V$ associated to
the eigenvalue $-1$. However, the transition map which sends
the hyperplane $m_2=m_3$ to the hyperplane $m_1=m_2$
is a composition of central projections (see the proof of Theorem~ \ref{thm:shapley}. In particular, it
sends the line through $V$ to the line through $V$ (remember
$V$ is an eigenvector) and it is a linear
fractional transformation on this line with derivative
$-1$ at the fixed point. But taking a map of the form
$f(x)=-x+\alpha x^2 +\gamma x^3+...$ we always
have $f^2(x)=x+0x^2+(?)x^3$. So in our case we get
that the second iterate of the linear fractional transformation
is equal to the identity map up to second order at the
fixed point. But this implies that it is in 
fact equal to the identity map.
So instead of a generic periodic doubling bifurcation, 
in this case we get non-generic dynamics at the bifurcation parameter
$\beta=\tau$; at this parameter value there is a one-parameter family 
(continuum) of periodic orbits. 

The computation showing that Shapley's periodic orbit 
is attracting for all $\beta\in (0,\sigma)$ goes similarly.
\end{proof}

\subsection{Dynamics and a periodic orbit on the set $J$}

We will now consider, again using $v^A, v^B$ space, 
the flow on the codimension two indifference set $J$ defined in the main
text. Recall that the flow which remains on $J$ is the unique
continuous extension to $J$ of the flow off $J$, everywhere except at the
interior equilibrium. 

\begin{prop}
The system has a Hopf-like bifurcation when $\beta=\sigma$
on the topological two-dimensional manifold $J$.
The orbits on $J$ spiral towards the interior 
equilibrium $E=(E^A,E^B)$ when $\beta\in (0,\sigma]$ while for
$\beta\in (\sigma,1)$ all orbits in $J$
(except the one remaining at $E$)
tend to a periodic orbit $\Gamma$ on $J$.
\end{prop}

\begin{figure}[htp]
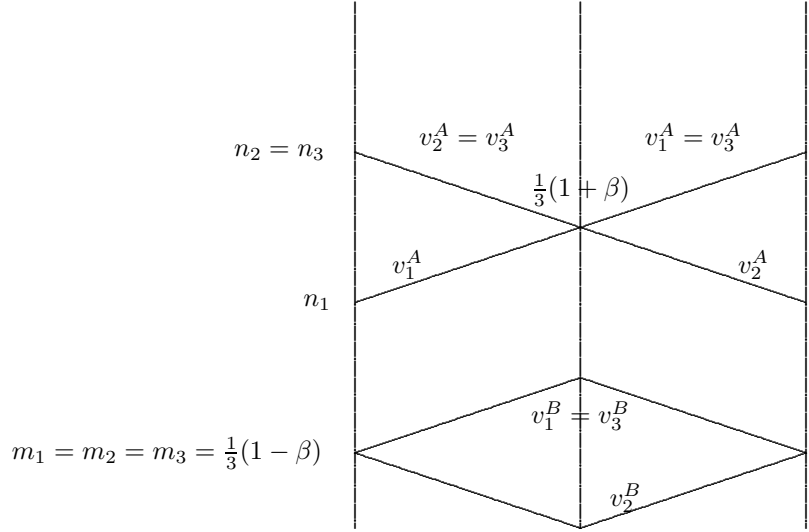
 \hfil
\beginpicture
\dimen0=0.1cm
\setcoordinatesystem units <\dimen0,\dimen0>
\setplotarea x from 0 to 60, y from 0 to 70
\setlinear
\plot 0 0 0 70 /
\plot 30 0 30 70 /
\plot 60 0 60 70 /
\put {$n_2=n_3$} at -10 50
\put {$n_1$} at -5 30
\put {$v^A_2=v^A_3$} at 15 52
\put {$v^A_1=v^A_3$} at 45 52
\put {$v^A_2$} at 53 35
\put {$v^A_1$} at 7 35
\put {$\frac{1}{3}(1+\beta)$} at 30 45
\plot 0 50 60 30 /
\plot 0 30 60 50 /
\put{$m_1=m_2=m_3=\frac{1}{3}(1-\beta)$} at -25 10
\plot 0 10 30 20 60 10 / \put {$v^B_1=v^B_3$} at 30 15
\plot 0 10 30 0 60 10 / \put {$v^B_2$} at 36 4
\endpicture
\caption{\label{figw}{\small Part of the
periodic orbit on $J$ with $(v^A, v^B)$ plotted against time.}}
\end{figure}

\begin{proof}
Again we compute a segment of orbit which makes one third of a complete
rotation around the interior equilibrium, and use the rotational
symmetry to compute the whole orbit.
We take an initial condition $v^A=n$ and $v^B=m$ where
$n_2=n_3 > n_1$ and 
$m_1=m_2=m_3=(1-\beta/3)$.  

Player A (whose preferences depend on $v^A$) is at first 
indifferent between strategies $2$ and $3$, and
so adjusts $p^A$ towards some convex combination of $P_2$ and $P_3$,
implying that $v^B$ moves towards the same convex combination of 
$B_2$ and $B_3$. The convex combination 
\begin{equation}
 \frac{(1+\beta)}{(2+\beta)} B_2 + \frac{1}{(2+\beta)} B_3 
=  \frac{1}{(2+\beta)}(1,\, -\beta(1+\beta),\, 1)
\label{eqn:app17}
\end{equation}
is chosen so that the orbit
remains in $J$, i.e.~so that $v^B_1(t)=v^B_3(t)$ in $t>0$. (We stress again
that this is not in any sense an arbitry choice; it is the unique
choice giving a continuous extension of the flow defined off the codimension
two indifference surface.) Consequently player $B$ will, for $t>0$, be
indifferent between strategies $1$ and $3$, will adjust $p^B$ towards
some convex combination of $P^B_1$ and $P^B_3$, and so $v^A$ will move
towards the same convex combination of $A_1$ and $A_3$. In this case, the
combination 
\begin{equation}
 \frac{1}{(1+\beta)} A_1 + \frac{\beta}{(1+\beta)} A_3 =
\frac{1}{(1+\beta)}(1+\beta^2,\beta,\beta)^{T} 
\label{eqn:app18}
\end{equation}
is chosen so that $v^A_2(t)=v^A_3(t)$ in $t>0$.

This motion continues until a time $t_1\in (0,1)$ at which
$v^A_1(t_1)=v^A_2(t_1)=v^A_3(t_1)=\frac{1}{3}(1+\beta)$
so that $t_1$ satisfies, using (\ref{eqn:app1}) and (\ref{eqn:app18})
$$v^A(t_1)=(1-t_1)v^A(0)+t_1\frac{1}{1+\beta}\left(\begin{array}{c}
1+\beta^2 \\ \beta \\ \beta \end{array} \right)=\frac{1+\beta}{3}
\left(\begin{array}{c} 1 \\ 1 \\ 1 \end{array} \right).$$
Hence
$$t_1=\frac{(n_2-n_1)(1+\beta)}{(n_2-n_1)(1+\beta)+(1+\beta^2-\beta)}\,.$$
We have also, from (\ref{eqn:app1}) and (\ref{eqn:app17}),    
$$v^B(t_1)=v^B(0)(1-t_1)+t_1 \frac{1}{2+\beta}
            \left( 1,\   -\beta-\beta^2,\  1 \right).$$
Now since $v^B_i(0)=(1-\beta)/3$ for $i=1,2,3$, we have
\begin{equation}
v^B_1(t_1)-v^B_2(t_1)=\frac{t_1}{2+\beta}[1+\beta+\beta^2]\,.
\label{eqn:app19}
\end{equation}

Player A now becomes indifferent between strategies $1$ and $3$,
while player B is still indifferent between strategies $1$ and $3$ 
(inspection shows this is the only way to remain on $J$). 
In the interval $t \in (t_1,t_1+t_2)$  we have that $v^A(t)$ heads
for 
\begin{equation}
   \frac{(1-\beta)}{(2-\beta)} A_1  + \frac{1}{(2-\beta)}A_3 =  
   \frac{1}{(2-\beta)} \left( 1,\, \beta(1-\beta),\, 1\right)^{T}
\label{eqn:app21}
\end{equation}
and $v^B(t)$ heads for
\begin{equation}
  \frac{(1+\beta)}{(1+2\beta)} B_1 + \frac{\beta}{(1+2\beta)}B_3 
 = \frac{1}{(1+2\beta)}\left(-\beta^2  ,\, 1+\beta,\, -\beta^2  \right) \,.        
\label{eqn:app22}
\end{equation}
Using (\ref{eqn:app1}) and (\ref{eqn:app22}) we seek to  
compute $t_2$ such that
$v^B(t_1+t_2)=\frac{(1-\beta)}{3}(1,1,1)$.
This gives
$$t_2 \left[ \frac{1}{(1+2\beta)}
(- \beta^2 ,\, 1+\beta ,\, -\beta^2 )-v^B(t_1) \right]=
\frac{1}{3}(1-\beta) (1,1,1) - v^B(t_1)$$
which implies
$$t_2=\frac{v^B_1(t_1)-v^B_2(t_1)}
{\frac{1+\beta+\beta^2}{1+2\beta} + v^B_1(t_1)-v^B_2(t_1)}=
\frac{t_1(1+2\beta)}{t_1(1+2\beta)+(2+\beta)}\,.$$
Now (\ref{eqn:app1}),  (\ref{eqn:app21}) and the fact that all components
of $v^A(t_1)$ are equal imply that
$$v^A_1(t_1+t_2)-v^A_2(t_1+t_2)
= \frac{t_2}{2-\beta}\left[1-\beta+\beta^2\right]\,.
$$
Substituting in for $t_2$ and then for $t_1$ we obtain
$$v^A_1(t_1+t_2)-v^A_2(t_1+t_2)
= \frac{t_1(1+2\beta)}{(t_1(1+2\beta)+(2+\beta))}
  \frac{(1-\beta+\beta^2)}{(2-\beta)}=$$
$$=\frac{(n_2-n_1)}
{(n_2-n_1)3(1+\beta)^2+ (2+\beta)(1-\beta+\beta^2)}
\frac{(1+\beta)(1+2\beta)(1-\beta+\beta^2)}{(2-\beta)}.$$
This formula gives us the value of $|v^A_1(t)-v^A_2(t)|$ at the end of 
the two legs, as a function $F()$ 
of the initial value $|n_2-n_1|$, where 
$$F(X):=\frac{X(1+\beta)(1+2\beta)(1-\beta+\beta^2)}
{(2-\beta)\left(3(1+\beta)^2 X + (2+\beta)(1+\beta^2-\beta)\right)}.$$
Now $F$ is a linear fractional transformation $F(X)=\frac{aX}{bX+c}$ and
$$F'(0)= \frac{a}{c}=\frac{(1+\beta)(1+2\beta)(1-\beta+\beta^2)}
{(2-\beta)(2+\beta)(1-\beta+\beta^2)}=
\frac{2\beta^2+3\beta+1}{-\beta^2+4}$$
which is less than $1$ when $\beta<\sigma$ and
greater than $1$ when $\beta>\sigma$. Since $F(X)\to 1$ as $X\to \infty$,
this implies that $x\mapsto F(x)$ has zero as an attracting fixed point
when $\beta\in (0,\sigma]$, and there exists an attracting fixed point
$x>0$ of $F$ when $\beta>\sigma$. This solution $x$
corresponds to the value of $n_2-n_1$ for which there is an
associated periodic orbit $\Gamma$ on $J$. This orbit attracts
all orbits on $J$ except the one starting
at the interior equilibrium $(E^A, E^B)$.

To prove the last past of the proposition, note that
$F(X)=X$ is equivalent
to
$$X=\frac{(1+\beta^2-\beta)(\beta^2+\beta-1)}
{(2-\beta)(1+\beta)^2}.$$
Now for the initial condition $n$ on $\Gamma$, 
we have $n_2-n_1=X=F(X)$, $n_2=n_3$, and  $n_1+n_2+n_3=(1+\beta)$.
Hence
$n_1=\frac{-2X+1+\beta}{3}$.
Note that
on the line $Z^A_{1,3}$ in $\Sigma_B$ (see Figure 9 in the main text)
one has (in $p^B$ coordinates)
$Q_{1,3}^B=\frac{1}{1+\beta}(\beta,1+\beta^2,\beta)$,
$E^B=\frac{1+\beta}{3}(1,1,1)$
and $R_{1,3}^B=\frac{1}{2-\beta}(1,\beta-\beta^2,1)$
(and symmetrically for other points).
During the third of the orbit studied  we have that we 
travel first along $[R^B_{23},E^B]$
until meeting $E^B$ and then along $[E^B,R_{13}^B]$.
The proportion
$$\frac{|E^B_{23}-f^B_{23}|}{|Q^B_{23}-f^B_{23}|}$$
is therefore, (by doing the computations in the $p^B$-simplex)
$$
\frac
{\frac{1+\beta}{3}-n_1}
{\frac{1+\beta^2}{1+\beta}-n_1}=
\frac
{\frac{1+\beta}{3}-\frac{-2X+1+\beta}{3}}
{\frac{1+\beta^2}{1+\beta}-\frac{-2X+1+\beta}{3}}=
\frac{\beta^2+\beta-1}{2\beta+1}.
$$
So when $\beta$ is equal to the golden mean
then this is equal to $0$
and it increases monotonically when $\beta\in [\sigma,1)$
(its maximum is $1/3$). Similarly,
$$\frac{|f^B_{23}-E^B_{23}|}{|R^B_{23}-E^B_{23}|}=
\frac
{\frac{1+\beta}{3}-n_1}
{\frac{1+\beta}{3}-\frac{\beta-\beta^2}{2-\beta}}=
\frac
{\frac{1+\beta}{3}-\frac{-2X+1+\beta}{3}}
{\frac{1+\beta}{3}-\frac{\beta-\beta^2}{2-\beta}}=
\frac{\beta^2+\beta-1}{(1+\beta)^2}$$
(which again increases in $\beta$).
\end{proof}

Note that when $\beta \ne \sigma$, the time taken for a trajectory
on $J$ to spiral into or out of the interior equilibrium $E$ is
finite. This is because for small $X$ we have that
$F(X) = c_1(\beta) X + O(X^2)$ and that the time taken for this
step is $t_1+t_2 = c_2(\beta) X + O(X^2)$ for some constants
$c_1(\beta),\, c_2(\beta)$ depending on $\beta$. When $c_1(\beta) \ne 1$,
which is true precisely when $\beta \ne \sigma$, this 
gives geometric (rather than exponential) convergence towards or 
divergence away from $E$. This phenomenon is intimately connected with the 
fact that the flow does not extend in a continuous way to $E$; there
is genuine non-uniqueness of the flow at $E$ (since trajectories may
remain at $E$ for all positive and negative time, or may leave it at
some moment in positive time ($\beta>\sigma$) or in negative 
time ($\beta<\sigma$).
When $\beta=\sigma$ the time taken for trajectories to tend to $E$ is
infinite, as can easily be checked.

\end{document}